\documentclass[12pt]{article} 
\newtheorem{proposition}{Proposition}[section] 
\usepackage{color} 
\usepackage{amssymb} 
\usepackage{amsmath} 
\usepackage{graphicx}
\usepackage{multirow} 
\usepackage{longtable} 
\usepackage{rotating} 
\usepackage{longtable} 
\usepackage{lscape} 
\usepackage[font=scriptsize]{caption} 
\date{} 
\title{Statistical Inference for Stable Distribution\\ Using EM algorithm\\
\vspace{1cm}
\small{{Mahdi Teimouri}}\\
Email: teimouri@aut.ac.ir\\
Department of Mathematics and Statistics, Faculty of Science and Engineering, \\
Gonbad Kavous University, Gonbad Kavous, Iran.} 
\begin{document}
\maketitle{} 
\noindent{} {\bf{Abstract:}} The class of $\alpha$-stable distributions with a wide range of applications in economics, telecommunications, biology, applied, and theoretical physics. This is due to the fact that it possesses both the skewness and heavy tails. Since $\alpha$-stable distribution suffers from a closed-form expression for density function, finding efficient estimators for its parameters has attracted a great deal of attention in the literature. Here, we propose some EM algorithm to estimate the maximum likelihood estimators of the parameters of $\alpha$-stable distribution. The performance of the proposed EM algorithm is demonstrated via comparison study in the presence of other well-known competitors and analyzing three sets of real data.\\\\
\noindent{} {\bf{Keyword:}} EM algorithm; Markov chain Monte Carlo; Maximum likelihood estimator; Price return; Profile likelihood; Stable distribution.
\section{Introduction}
Despite the lack of non-analytical expression for the probability density function (pdf), the class of $\alpha$-stable distributions are becoming increasingly popular in such fields as economics, finance, and insurance. Details for the applications of $\alpha$-stable distributions in these fields can be found in \cite{Mandelbrot2006}, \cite{Mittnik2003}, \cite{Nolan2003}, \cite{Ortobelli2010}, \cite{Rachev2003},  \cite{Rachev2005}, and \cite{Rachev2000}. The only exceptions for which the pdf has analytical expression are Gaussian, Levy, and Cauchy distributions. This can be regarded as a major obstacle in the way of using this class in practice. This is while the characteristic function (chf) of $\alpha$-stable distributions has closed-form expression and takes different forms, see \cite{Nolan1998} and \cite{Zolotarev1986}. In what follows we review briefly two important forms which are known in the literature as $S_{0}$ and $S_{1}$ parameterizations. If random variable $Y$ follows $\alpha$-stable distribution, then the chf of $Y$, i.e., $\varphi_{Y}{(t)}=E \exp (jtY)$, in $S_{0}$ parameterization is given by
\begin{eqnarray} 
\label{chf0} 
\varphi_{Y}{(t)}=\left\{\begin{array}{*{20}c} 
\exp\biggl\{-\left| \sigma t \right|^\alpha\Big[1-j\beta~\mathrm{sgn}(t)\tan\left(\frac{\pi \alpha}{2}\right)\Big]+jt\mu_{0} -jt\beta \sigma \tan\left(\frac{\pi \alpha}{2}\right)\biggr\},~ \mathrm{{if}}~\alpha \ne 1,\\ 
\exp\biggl\{-\left| \sigma t \right| \Big[1+j\beta~\mathrm{sgn}(t)\frac{2}{\pi}\log \left| t \right|\Big]+j t\mu_{0}-\frac{2}{\pi}t\beta \sigma \log \sigma \biggr\},~~~~~~~~~\mathrm{if}~\alpha= 1. 
\end{array} \right. 
\end{eqnarray} 
where $j^{2}$=-1 and $\mathrm{sgn}(.)$ is the well-known sign function. The family of stable distributions has four parameters: tail thickness $\alpha \in (0, 2]$, skewness $\beta \in [-1, 1]$, scale $\sigma \in \mathbb{R}^{+}$, and location $\mu_{0} \in \mathbb{R}$. If $\beta$=0, it would be the class of symmetric $\alpha$-stable (S$\alpha$S) distributions. If $\beta$=1 and $\alpha< 1$, we have the class of the positive stable distributions. In this case, $Y$ varies over the positive semi-axis of real line. Also, the chf of $Y$ in $S_{1}$ parameterization is given by
\begin{eqnarray} 
\label{chf1} 
\varphi_{Y}{(t)}=\left\{\begin{array}{*{20}c} 
\exp\biggl\{-\left| \sigma t \right|^\alpha\Big[1-j\beta~\mathrm{sgn}(t)\tan\left(\frac{\pi \alpha}{2}\right)\Big]+jt\mu_{1} \biggr\},~~~~\mathrm{{if}}~\alpha \ne 1,\\ 
\exp\biggl\{-\left| \sigma t \right| \Big[1+j\beta~\mathrm{sgn}(t)\frac{2}{\pi}\log \left| t \right|\Big]+j t\mu_{1} \biggr\},~~~~~~~\mathrm{if}~\alpha= 1, 
\end{array} \right. 
\end{eqnarray} 
where $\mu_{1} \in \rm I\!R$ is the location parameter. The pdf and chf in $S_{0}$ parameterization are continuous functions of all four parameters, while the pdf and chf in $S_{1}$ parameterization both has discontinuity at $\alpha=1$. As it is seen, the parameterizations (\ref{chf1}) and (\ref{chf0}) differs only for the location parameter. The location parameters in $S_{0}$ and $S_{1}$ parameterizations, respectively shown by $\mu_{0}$ and $\mu_{1}$, are related as
\begin{eqnarray} 
\label{chf} 
\mu_{1}=\left\{\begin{array}{*{20}c} 
\mu_{0}-\beta \sigma \tan\left(\frac{\pi \alpha}{2}\right),~~~~~\mathrm{{if}}~\alpha \neq1,\\ 
\mu_{0}-\beta\frac{2}{\pi} \sigma \log \sigma,~~~~~~~~~\mathrm{if}~ \alpha= 1.\\ 
\end{array} \right. 
\end{eqnarray} 
Clearly, when $\beta=0$ both $S_{0}$ and $S_{1}$ parameterizations are equal. \par Hereafter, we write $S_{0}(\alpha,\beta,\sigma,\mu_{0})$ and $S_{1}(\alpha,\beta,\sigma,\mu_{1})$ to denote the class of stable distributions in $S_{0}$ and $S_{1}$ parameterizations, respectively. For the class of normal distributions with mean $a$ and variance $b$, i.e., ${\cal{N}}(a,b)$, the pdf at point $y$ is shown by $\phi(y,a,b)$. The generic symbol ${\cal{E}}(\lambda)$ accounts for the family of exponential distributions with rate parameter $\lambda>0$, and ${\cal{W}}(a,b)$ denotes a Weibull distribution with pdf $ab^{-a}w^{a-1}\exp\bigl\{-({w/b})^{a}\bigr\}$; for $w>0$, shape parameter $a>0$, and scale parameter $b>0$. Also $f_{Y}(.|\alpha,\beta,\sigma,\mu_{0})$, $g(.|\alpha)$, and $h(.|\alpha)$ are pdfs of distributions $S_{0}(\alpha,\beta,\sigma,\mu_{0})$, $S_{1}\bigl(\frac{\alpha}{2},1,\bigl(\cos(\pi \alpha/4)\bigr)^{2/\alpha},0\bigr)$, and $S_{1}(\alpha,1,1,0)$, respectively.
\par
 This paper is organized as follows. In what follow firstly the EM algorithm and its extensions are reviewd briefly, then some properties of $\alpha$-stable distributions are given. In Section 2, the proposed EM algorithm is introduced. Model validation using simulations and real data analysis is carried out in Section 3. Some conclusions are made in Section 4.
\subsection{EM algorithm and its extensions}
The EM algorithm is the most popular approach for estimating the parameters of a statistical model when we encounter missing (or latent) observations, see \cite{Dempster1977}. Assume that $\underline{\boldsymbol{x}}=({\boldsymbol{x}}_{1},\dots,{\boldsymbol{x}}_{n})$ denotes the vector of complete data with pdf $d(\underline{\boldsymbol{x}}|\Theta)$ in which $\boldsymbol{x}_i=(y_i,\xi_i)$ consists of observed and missing values, respectively. Also, $L_c(\Theta)=\Pi_{i=1}^{n}d(\underline{{x}_i}|\Theta)$ denotes the likelihood function of complete data in which $\Theta$ is the parameter vector. The EM algorithm finds the parameter vector that maximizes the conditional expectation of the log-likelihood of complete data given the observed data and a current guess $\Theta^{(t)}$, of the parameter vector. Usually, the EM algorithm works as follows.
\begin{enumerate}
\item E-step: given $\underline{y}=(y_1,\dots,y_n)$ and $\Theta^{(t)}$, it computes $Q\bigl(\Theta \big|\Theta^{(t)}\bigr)=E\bigl(l_c(\Theta)\big|\underline{y}, \Theta^{(t)}\bigr)$.
\item M-step: it finds such $\Theta$ that $Q\bigl(\Theta \big|\Theta^{(t)}\bigr)$ is maximized.
\end{enumerate}
Notice that $l_c(\Theta)$ refers to the logarithm of $L_c(\Theta)$. Both steps of the EM algorithm are repeated until convergence occurs, see \cite{McLachlan2008}. The convergence of the EM algorithm is guaranteed by \cite{Little}. In what follows we describe two extensions of the EM algorithm.
\subsubsection{ECM Algorithm} 
When implementing the M-step of the EM algorithm is mathematically impossible, an additional step is considered. Besides the E- and M-step, a sequence of the conditional maximization is carried out. The new step is known as the CM-step and the algorithm is known as the ECM algorithm. In the ECM algorithm, the CM-step finds the maximum of $Q\bigl(\Theta \big|\Theta^{(t)}\bigr)$ through maximizing the constrained marginal log-likelihood function, see \cite{Liu1994} and \cite{Meng1993}.
\subsubsection{Stochastic EM Algorithm} 
For a complete data of size $n$, assuming we are currently at $t$th iteration, each stochastic EM (SEM) algorithm amounts to a four-step sequence given by the following.
\begin{enumerate}
\item Given $\Theta^{(t)}$ and $y_{i}$, a sequence of latent (or missing) variables, i.e., $\underline{\xi}=(\xi_{1},\dots,\xi_{n})$ is simulated from posterior pdf $d\bigl(\xi_i\big|\Theta^{(t)}, y_{i}\bigr)$; for $i=1,\dots, n$. 
\item The simulated latent realizations are replaced into the log-likelihood function of complete data. 
\item The EM algorithm is applied to the set of complete data, $\underline{{\boldsymbol{x}}}=({\boldsymbol{x}}_{1},\dots,{\boldsymbol{x}}_{n})$ in which ${\boldsymbol{x}_i}=(y_i,\xi_i)^{T}$; for $i=1,\dots,n$, is the vector of observed and latent realizations. 
\item The vector of parameters are updated as $\Theta^{(t+1)}$ and then it is used to simulate from posterior pdf $d\bigl(\xi_{i}\big|\Theta^{(t+1)}, y_{i}\bigr)$; for $i=1,\dots,n$.
\end{enumerate}
Under some mild regularity conditions, by repeating above four-step process for a sufficiently large number of cycles, the distribution of $\bigl\{\Theta^{(t+1)}\bigr\}$ constitutes a Markov chain that converges to a stationary distribution, see \cite{Broniatowski1983}, \cite{Celeux1985}, and \cite{Ip1994}. Contrary to the EM algorithm, for the SEM algorithm convergence occurs in distribution and in practice the number of cycles is determined through a graphical display. Suppose $M_0$ to be the length of burn-in period and $M$ is a sufficiently large number, the SEM algorithm estimates $\Theta$ as:
\begin{eqnarray} 
\hat{\Theta}=\frac{1}{M-M_0}\sum_{t=M_0+1}^{M} \Theta^{(t)}.\nonumber
\end{eqnarray}
\subsection{Some properties of $\alpha$-stable distribution} 
Estimating the parameters of a $\alpha$-stable distribution through the EM algorithm needs a hierarchy or stochastic representation. Here we give three useful results which play a large role to follow the next chapter.
\begin{proposition}\label{prop1} Suppose $Y\sim S_{0}(\alpha,\beta,\sigma,\mu_{0})$, $P \sim S_{1}\bigl(\alpha/2,1,\bigl(\cos(\pi \alpha/4)\bigr)^{2/\alpha},0\bigr)$, and $V\sim S_{1}(\alpha,1,1,0)$. We have
\begin{eqnarray} \label{asymm}
Y\mathop=\limits^d\eta\sqrt{2P}N+ \theta V+\mu_{0}- \lambda
\end{eqnarray}
where $\mathop=\limits^d$ denotes the equality in distribution, $\eta=\sigma\left(1-|\beta|\right)^{\frac{1}{\alpha}}$, $\theta=\sigma\mathrm{sgn}(\beta)|\beta|^{\frac{1}{\alpha}}$, $\lambda=\sigma \beta \tan \bigl(\pi \alpha/2 \bigr)$, and $N\sim{\cal{N}}(0,1)$. All random variables $N$, $P$, and $V$ are mutually independent.
\end{proposition}
\begin{proposition} \label{prop2} Let $Y\sim S_{0}\bigl(\alpha,\beta,\sigma,\mu_{0})$ be independent of $V\sim S_{1}\bigl(\alpha,1,1,0)$. Then, 
\begin{align*} 
\frac{Y-\theta V-\mu_{0}+\lambda}{\delta} \sim S_{1}(\alpha,0,1,0)
\end{align*} 
where $\theta=\sigma\mathrm{sgn}(\beta)|\beta|^{\frac{1}{\alpha}}$, $\delta=\sigma(1+ |\beta|)^{\frac{1}{\alpha}}$, and $\lambda=\sigma \beta \tan \bigl(\pi \alpha/2 \bigr)$.
\end{proposition}
\begin{proposition} \label{prop3} Let $E\sim {\cal{E}}(1)$ be independent of $S\sim S_{1}\bigl(\alpha,0,1,0)$. Then, 
\begin{align*} 
\frac{S}{\sqrt{2E}} \mathop=\limits^d\frac{N}{W}.
\end{align*} 
where $N \sim {\cal{N}}(0,1)$ and $W \sim {\cal{W}}(\alpha,1)$.
\end{proposition}
\section{Proposed EM algorithm}
Assume that $y_{1},\dots, y_{n}$ constitute a sequence of identically and independent realizations of $\alpha$-stable distribution in $S_{0}$ parameterization. The vector of complete data associated with (\ref{asymm}) is shown by $\underline{\boldsymbol{x}}=(\underline{x}_1,\dots,\underline{x}_n)=\bigl((y_1,p_1,{v}_1),\dots,(y_n,p_n,{v}_n)\bigr)$ in which $\underline{p}$, $\underline{v}$ are vectors of realizations of latent variables $\underline{P}$ and $\underline{V}$, respectively. It turns out that representation (\ref{asymm}) admits the following hierarchy. 
\begin{align}\label{mixhier}
Y|P=p,V=v\sim& {\cal{N}}\bigl(\mu_{0}-\lambda+\theta v,2p\eta^{2}\bigr),\nonumber\\
P\sim& S_{1}\Bigl(\frac{\alpha}{2},1,\bigl(\cos(\pi \alpha/4)\bigr)^{\frac{2}{\alpha}},0\Bigr),\nonumber\\
V\sim& S_{1}\bigl(\alpha,1,1,0\bigr),
\end{align}
where $\eta$, $\theta$, and $\lambda$ are defined after Proposition \ref{prop1}. Using hierarchy (\ref{mixhier}), the log-likelihood function of complete data is
\begin{align}
l_{c}(\Theta)=&\text{C}-n\log \eta-\frac{1}{2}\sum_{i=1}^{n}\biggl(
\frac{y_{i}-\mu_{0}+\lambda-\theta v_{i}}{\sqrt{2p_{i}}\eta}\biggr)^{2}+\sum_{i=1}^{n}\log f_{P_{i}}(p_i|\alpha)+\sum_{i=1}^{n}\log f_{V_{i}}(v_i|\alpha),\nonumber
\end{align}
where $\text{C}$ is a constant independent of the parameters vector $\Theta=(\alpha,\beta,\sigma,\mu_{0})^{T}$. The conditional expectation of $l_{c}(\Theta)$, i.e., $Q\bigl(\Theta \big|\Theta^{(t)}\bigl)=E\bigl(l_c(\Theta;\boldsymbol{p}, \boldsymbol{v})\big|\boldsymbol{y}, \Theta^{(t)}\bigr)$ is
\begin{align}
Q\bigl(\Theta\big|\Theta^{(t)}\bigr)=&\text{C}-n\log \eta-\frac{\theta^{2}}{4\eta^{2}}\sum_{i=1}^{n}
E^{(t)}_{2i}+\frac{\theta}{2\eta^{2}}\sum_{i=1}^{n}\bigl(y_{i}-\mu_{0}+\lambda \bigr)
E^{(t)}_{1i}\nonumber\\
&-\frac{1}{4\eta^{2}}\sum_{i=1}^{n}\bigl(y_{i}-\mu_{0}+\lambda\bigr)^{2}
E^{(t)}_{0i}+\sum_{i=1}^{n}E\bigl(\log f_{P_{i}}(p_i|\alpha)\bigr)+\sum_{i=1}^{n}E\bigl(\log f_{V_{i}}(v_i|\alpha)\bigr).\nonumber
\end{align}
To complete the E-step, we need to compute the quantities $E^{(t)}_{ri}=E\bigl(P_{i}^{-1}V_{i}^{r}\big| y_{i}, \Theta^{(t)}\bigr)$; for $r=0,1,2$. Now, it is straightforward to show that
\begin{align} 
\label{expectation} 
E^{(t)}_{ri}=& 
\frac{1}{2\eta^{(t)}\sqrt{\pi} f_{Y}\bigl(y_i\big |\alpha^{(t)},{\beta}^{(t)}, {\sigma}^{(t)}, \mu_{0}^{(t)}\bigr)}\nonumber\\
&\times \int \int p^{-1.5}v^{r}\exp\biggl\{-\frac{1}{2}\Bigl(\frac{y_{i}-\mu^{(t)}_{0}+\lambda^{(t)}-\theta^{(t)} v}{\sqrt{2p}\eta}\Bigr)^{2}\bigg\}
h\bigl(v|\alpha^{(t)}\bigr)g\bigl(p|\alpha^{(t)}\bigr)dvdp,
\end{align} 
where $\eta^{(t)}=\sigma^{(t)}\left(1-|\beta^{(t)}|\right)^{\frac{1}{\alpha^{(t)}}}$, $\theta^{(t)}=\sigma^{(t)} \mathrm{sgn}(\beta^{(t)})|\beta^{(t)}|^{\frac{1}{\alpha^{(t)}}}$, and $\lambda^{(t)}=\beta^{(t)} \sigma^{(t)} \tan\bigl(\frac{\pi \alpha^{(t)}}{2}\bigr)$. Details for evaluating $E^{(t)}_{ri}$ in (\ref{expectation}) is described in Appendix 1. Assuming we are currently at $t$th iteration, all steps of the proposed EM algorithm including the E-, M-, CM-step, and maximizing the profile log-likelihood function are given by the following.
\begin{itemize}
\item {\bf{E-step}}: Given current guess of $\Theta$, i.e., $\Theta^{(t)}$, the quantity $E^{(t)}_{ri}$; for $i=1,\dots,n$ is computed.
\item {\bf{M-step}}: Given $\Theta^{(t)}$, the parameter vector $\Theta$ is updated as $\Theta^{(t+1)}$ by maximizing $Q\bigl(\Theta\big|\Theta^{(t)}\bigr)$ with respect to $\Theta$. In the M-step, the location parameter is updated as 
\begin{align}
\mu_{0}^{(t+1)}=\frac{\sum_{i=1}^{n}\bigl(y_{i}+\lambda^{(t)}\bigr)E^{(t)}_{0i}+\eta^{(t)}\sum_{i=1}^{n}E^{(t)}_{1i}}{\sum_{i=1}^{n}E^{(t)}_{0i}},\nonumber
\end{align}
and the updated scale parameter $\sigma^{(t+1)}$, is the root of equation ${\cal{G}}(\sigma)=a\sigma^2+b\sigma+c$ in which $a=-n$,
\begin{align}
b=&\frac{\beta^{(t)} \tan\bigl(\frac{\pi \alpha^{(t)}}{2}\bigr)\sum_{i=1}^{n}\bigl(y_{i}-\mu_{0}^{(t+1)}+\lambda^{(t)}\bigr)E^{(t)}_{0i}}{2\bigl(1-|\beta^{(t)}|\bigr)^{\frac{2}{\alpha^{(t)}}}}\nonumber\\
&-\frac{\mathrm{sgn}(\beta^{(t)})|\beta^{(t)}|^{\frac{1}{\alpha^{(t)}}}\sum_{i=1}^{n}\bigl(y_{i}-\mu_{0}^{(t+1)}+\lambda^{(t)}\bigr)E^{(t)}_{1i}}{2\bigl(1-|\beta^{(t)}|\bigr)^{\frac{2}{\alpha^{(t)}}}},\nonumber
\end{align}
and
\begin{align}
c=\frac{\sum_{i=1}^{n}\bigl(y_{i}-\mu_{0}^{(t+1)}+\lambda^{(t)}\bigr)^{2}E^{(t)}_{0i}}{2\bigl(1-|\beta^{(t)}|\bigr)^{\frac{2}{\alpha^{(t)}}}}.\nonumber
\end{align}
The updated skewness parameter is obtained as $\beta^{(t+1)}=\arg\max\limits_{\beta} ~{\cal{F}} (\beta)$ where 
\begin{align}\label{betaupdating}
{\cal{F}} (\beta)=&-\frac{n\log(1-|\beta|)}{\alpha^{(t)}}+
\frac{1}{4}\Bigl(\frac{|\beta|}{1-|\beta|}\Bigr)^{\frac{2}{\alpha^{(t)}}}
\sum_{i=1}^{n}E^{(t)}_{2i}\nonumber\\
&+\frac{\mathrm{sgn}(\beta)|\beta|^{\frac{1}{\alpha^{(t)}}}}{2\sigma^{(t+1)}
\bigl(1-|\beta|\bigr)^{\frac{2}{\alpha}}}
\sum_{i=1}^{n}\Bigl(y_{i}-\mu_{0}^{(t+1)}+ \beta \sigma^{(t+1)} \tan\bigl(\pi \alpha^{(t)}/2\bigr) \Bigr)E^{(t)}_{1i}\nonumber\\
&+\frac{1}{4\bigl(1-|\beta|\bigr)^{\frac{2}{\alpha}}}
\sum_{i=1}^{n}\Bigl(y_{i}-\mu_{0}^{(t+1)}+ \beta \sigma^{(t+1)} \tan\bigl(\pi \alpha^{(t)}/2\bigr) \Bigr)^{2}E^{(t)}_{0i}.
\end{align}
\item {\bf{CM-step}}: 
The tail thickness parameter is updated in the CM-step by maximizing the marginal log-likelihood function with respect to $\alpha$ as 
\begin{align*} 
\alpha^{(t+1)}=\arg\max\limits_{\alpha} ~\sum_{i=1}^{n}\log f\bigl(y_i\big|\alpha,{\beta}^{(t+1)}, {\sigma}^{(t+1)}, \mu_{0}^{(t+1)}\bigr)
\end{align*} 
For this, at $t$-th iteration of the EM algorithm, we apply the SEM algorithm by three following steps to obtain $\alpha^{(t+1)}$.
\begin{enumerate} 
\item Let $e_1, e_2,\dots,e_n$ be independent and identically distributed realizations from ${\cal{E}}(1)$. Consider the transformation ${{y}}_{i}^{**}={{y}}_{i}^{*}/\sqrt{2e_i}$ in which 
\begin{align*} 
y^{*}_{i}=\frac{y_{i}-\theta^{(t+1)} v_{i}-{\mu}^{(t+1)}_{0}-\lambda^{(t+1)}}{\delta^{(t+1)}}
\end{align*} 
where $\delta^{(t+1)}=\sigma^{(t+1)}\bigl(1+ |\beta^{(t+1)}|\bigr)^{\frac{1}{\alpha^{(t)}}}$ and $\lambda^{(t+1)}=\beta^{(t+1)} \sigma^{(t+1)} \tan\bigl(\frac{\pi \alpha^{(t)}}{2}\bigr)$; for $i=1,\dots,n$. It turns out, from Propositions \ref{prop2} and \ref{prop3}, that 
\begin{align} \label{hierarchy2}
Y^{**}_i\big|W_i=w_i &\sim {\cal{N}}\bigl(0, w^{-2}_{i}\bigr),\nonumber\\ 
W_i &\sim {\cal{W}}(\alpha,1).
\end{align}
Based on hierarchy (\ref{hierarchy2}), the log-likelihood of complete data can be written as 
\begin{align}\label{lw} 
l_{c}\bigl(\alpha\big|y_i^{**}\bigr)=&\text{C}+n \log \alpha-\sum_{i=1}^{n}w^{\alpha}_i
+\alpha \sum_{i=1}^{n}\log w_i. 
\end{align}
\item Considering $W$ as the latent variable, we simulate $\underline{w}=(w_{1},\dots,w_{n})$ from posterior distribution $W_{i}$ given $Y^{**}_i$ and $\alpha^{(t)}$; for $i=1,\dots,n$, using the method described in Appendix 3.
\item Substitute $\underline{w}$ in right-hand side of (\ref{lw}) and maximize it with respect to $\alpha$ to obtain $\alpha^{(t+1)}$. Obtaining $\alpha^{(t+1)}$, we go back to the step one and repeat the CM-step for $M$ cycles. This yields a sequence of $M$ updated tail thickness parameter as: $\alpha^{(t+1,1)},\dots,\alpha^{(t+1,M)}$. Now, the tail thickness parameter is updated as
\begin{align*} 
\alpha^{(t+1)}=\frac{1}{M-M_0}\sum_{j=M_0+1}^{M} \alpha^{(t+1,j)},
\end{align*} 
where $M_0$ in the length of burn-in period of the SEM and $\alpha^{(t,j)}$ is updated tail thickness parameter at $j$th cycle of the CM-step while we are at $t$-th iteration of the EM algorithm. Once we obtain $\alpha^{(t+1)}$, we go back and perform the EM algorithm from E-step for a sufficiently large number of iterations, say $N$. After a burn-in period of length $N_0$, the EM algorithm converges to the true distribution. 
\end{enumerate}
Updating the skewness parameter by maximizing ${\cal{F}}(\beta)$ in (\ref{betaupdating}) yields result that goes to zero. So, we update $\beta^{(t)}$ using optimization tools such as $\mathsf{optim}$ developed in $\mathsf{R}$ package by maximizing profile log-likelihood function. For this, after obtaining EM-based estimations of ${\alpha}$, ${\sigma}$, and ${\mu_{0}}$, namely, $\hat{\alpha}_{EM}$, $\hat{\sigma}_{EM}$, and $\hat{\mu_{0}}_{EM}$, we maximize $\sum_{i=1}^{n}\log f_{Y} \bigl(y_i\big|\hat{\alpha}_{EM}, \beta, \hat{\sigma}_{EM}, \hat{\mu_{0}}_{EM}\bigr)$ with respect to $\beta$ to obtain $\hat{\beta}_{EM}$. Fortunately, this approach leads to satisfactory results.
\end{itemize}
\section{Model validation using simulated and real data}
Here, firstly, we perform a simulation study to compare the performance of the proposed EM algorithm with the ML approach for estimating the parameters of the class $S_{0}(\alpha,\beta,\sigma,\mu_{0})$. To do this, data are generated by the method of simulating $\alpha$-stable random variable (see \cite{Chambers1976}) and then the maximum likelihood (ML) estimations are evaluated using $\mathsf{STABLE}$ software, see \cite{Nolan2001}. Secondly, we apply the proposed EM and ML approaches to the five sets of real data.
\subsection{Simulation study}
\begin{figure}
\resizebox{\textwidth}{!}
{\begin{tabular}{ccc}
\includegraphics[width=40mm,height=40mm]{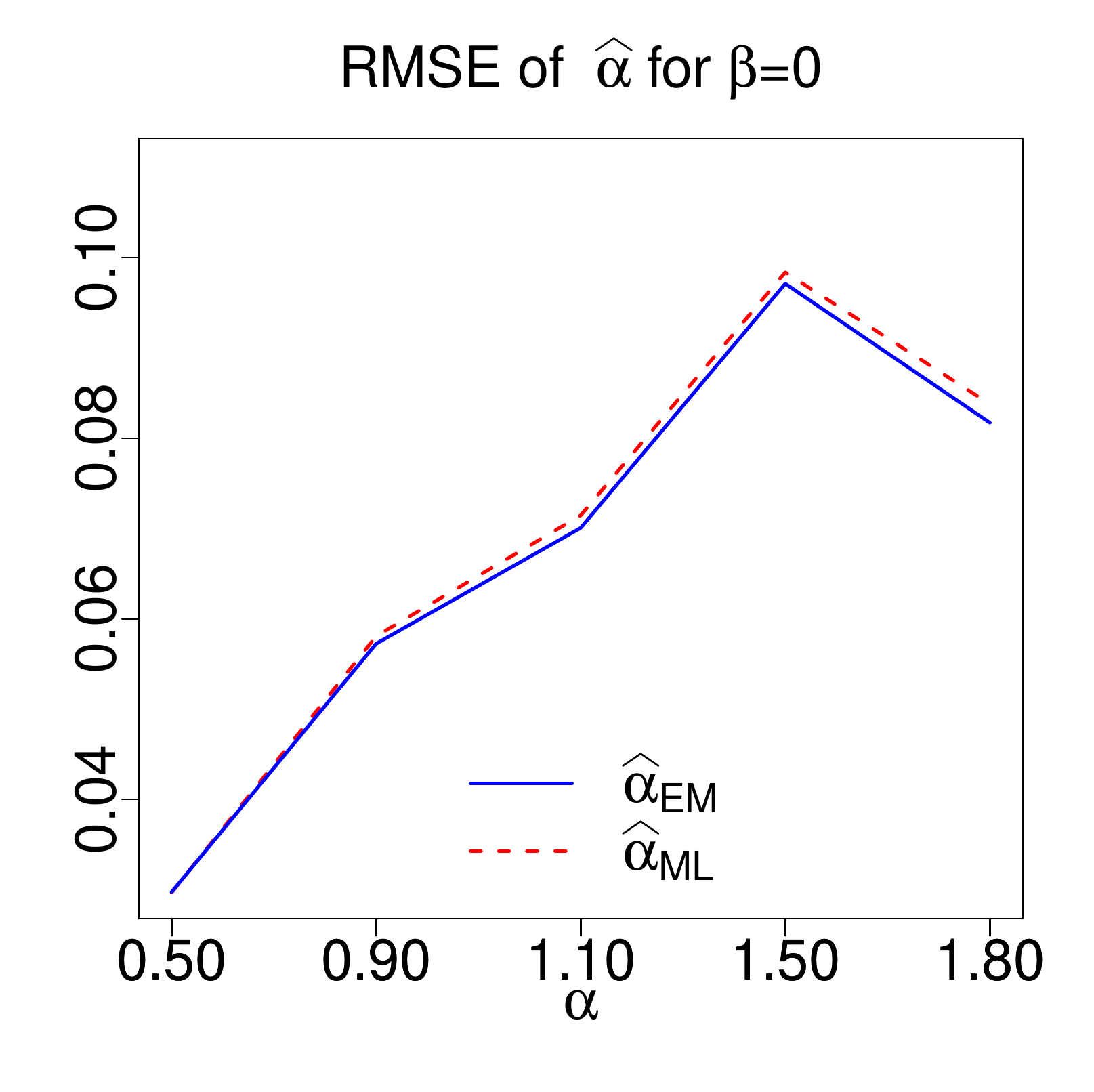}&
\includegraphics[width=40mm,height=40mm]{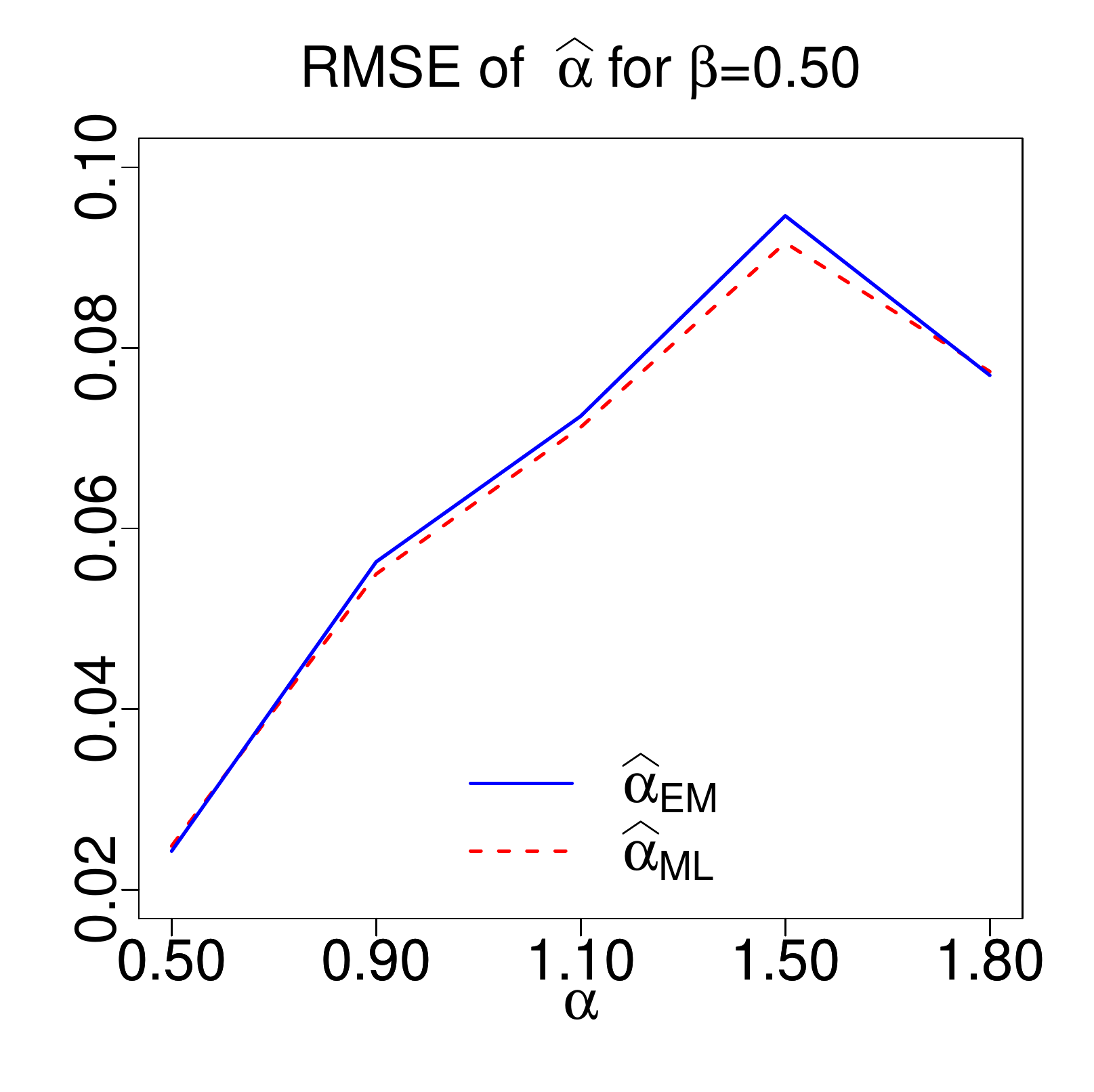}&
\includegraphics[width=40mm,height=40mm]{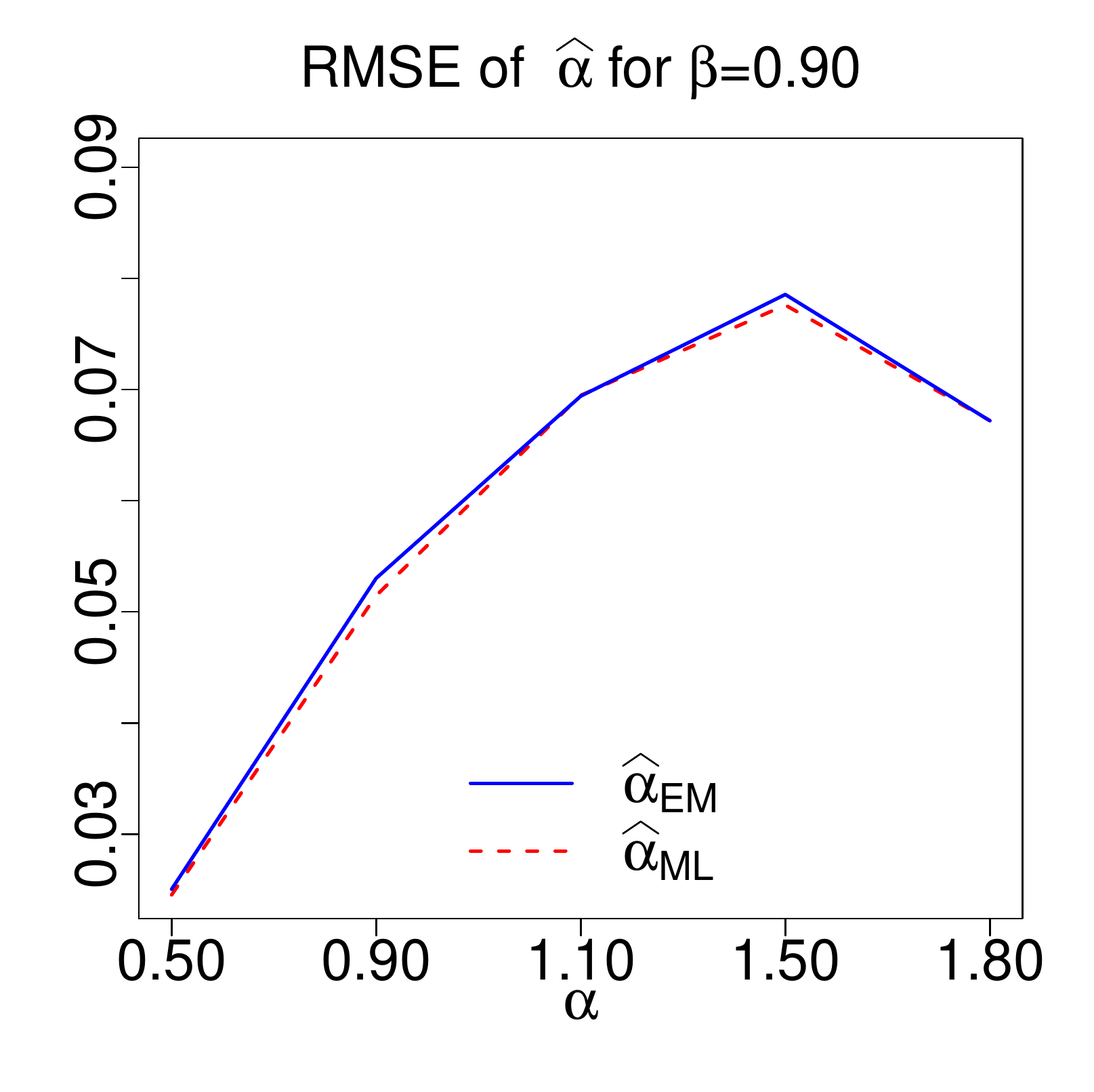}\\
\includegraphics[width=40mm,height=40mm]{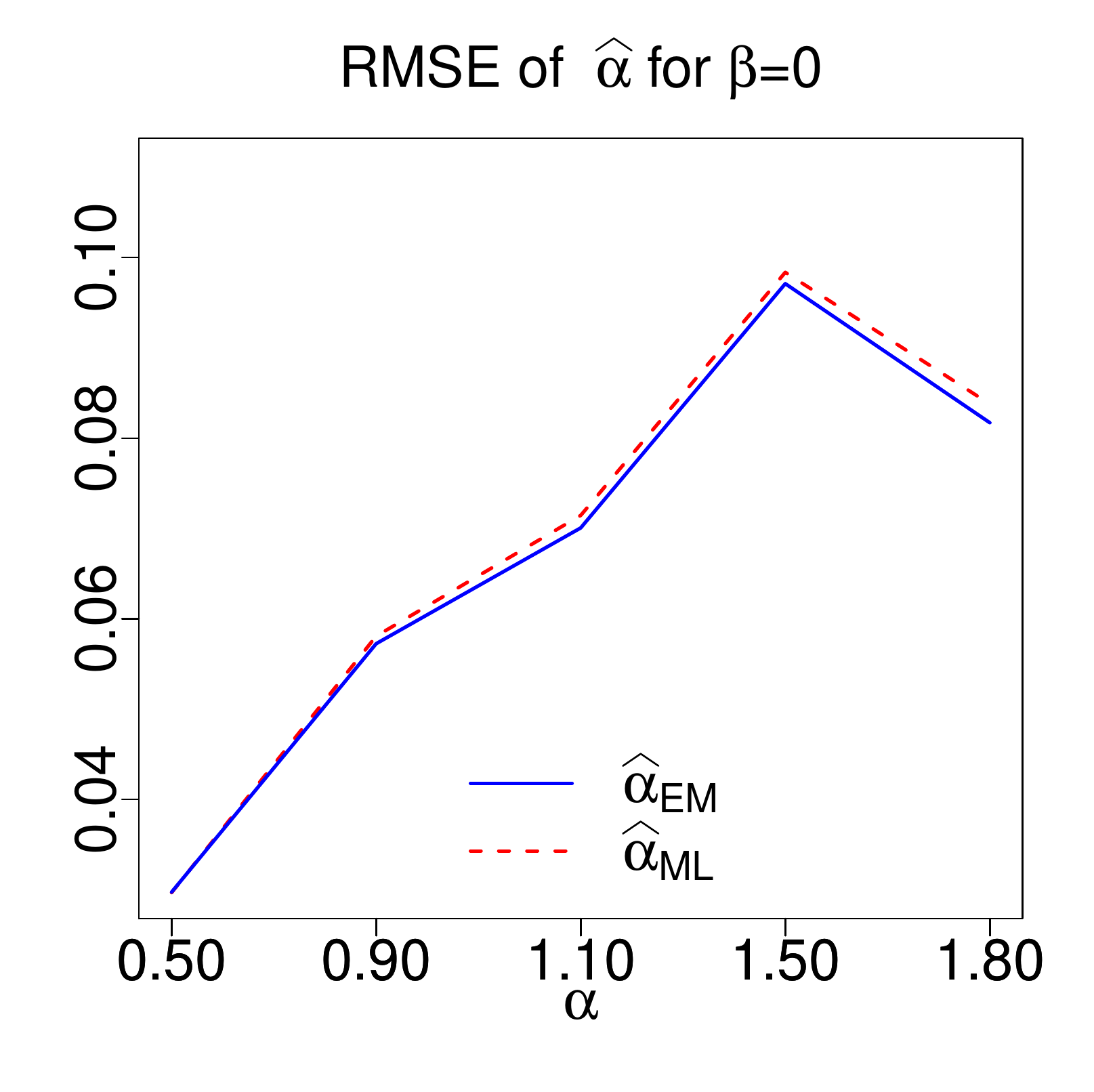}&
\includegraphics[width=40mm,height=40mm]{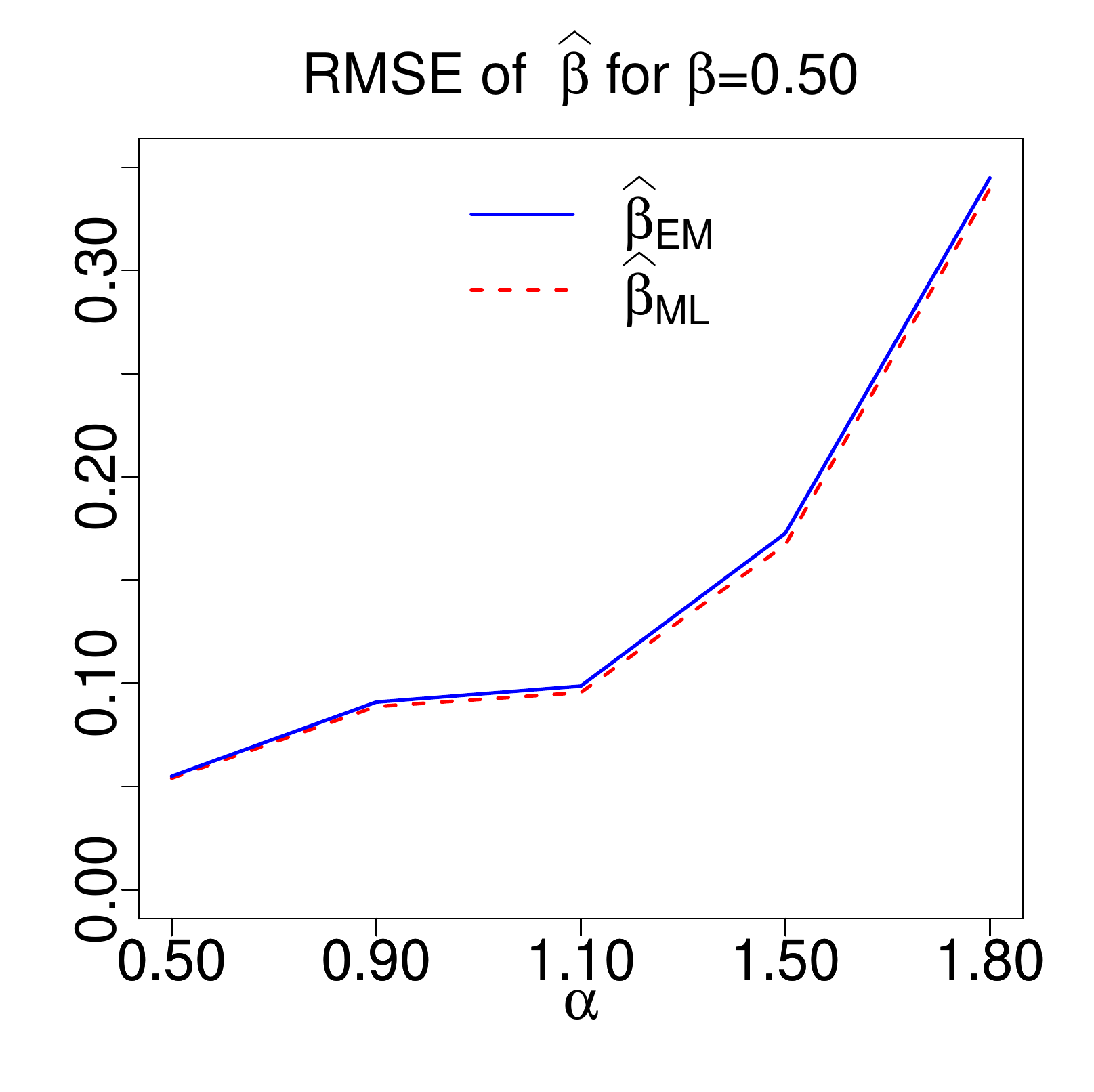}&
\includegraphics[width=40mm,height=40mm]{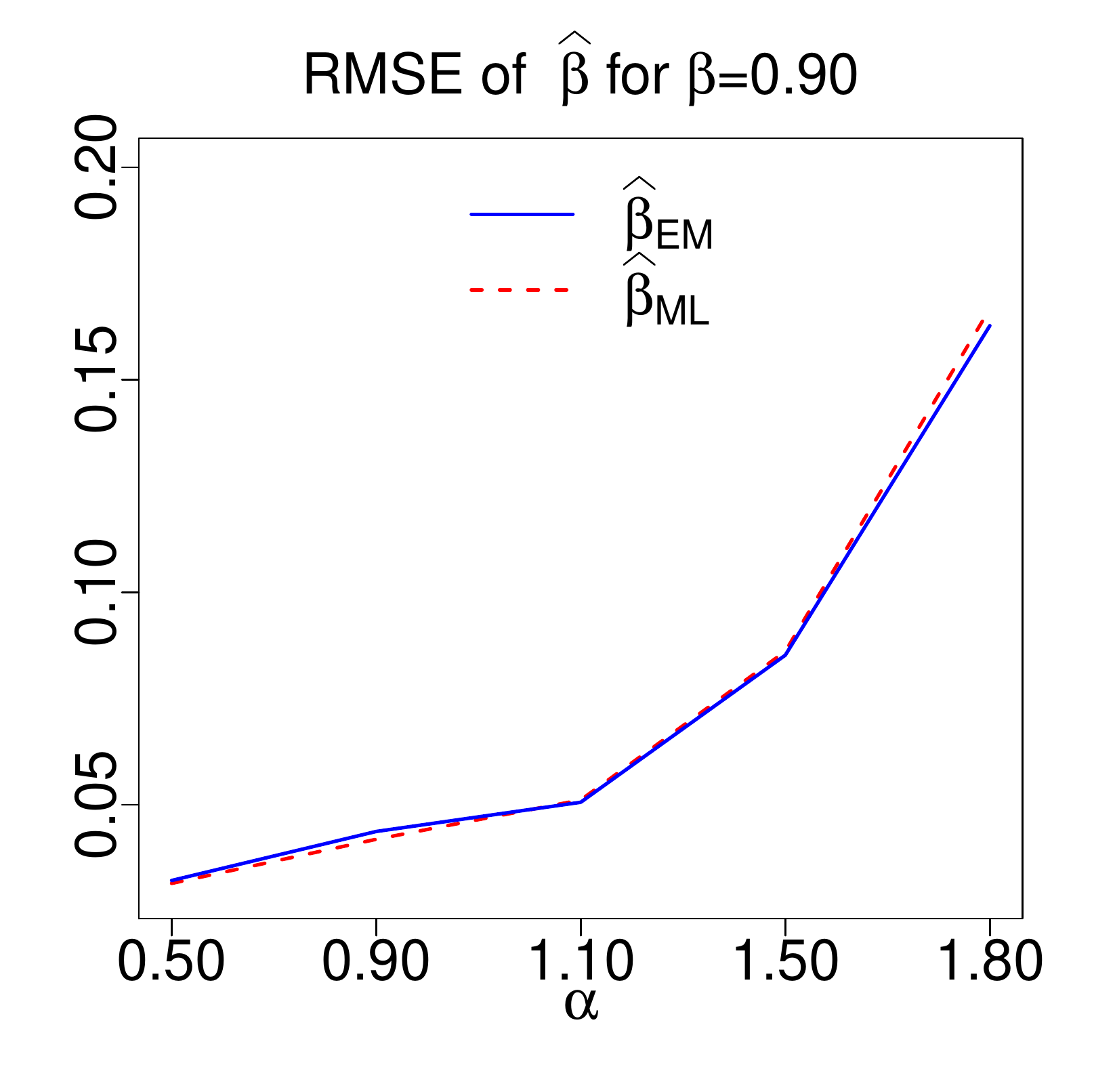}\\
\includegraphics[width=40mm,height=40mm]{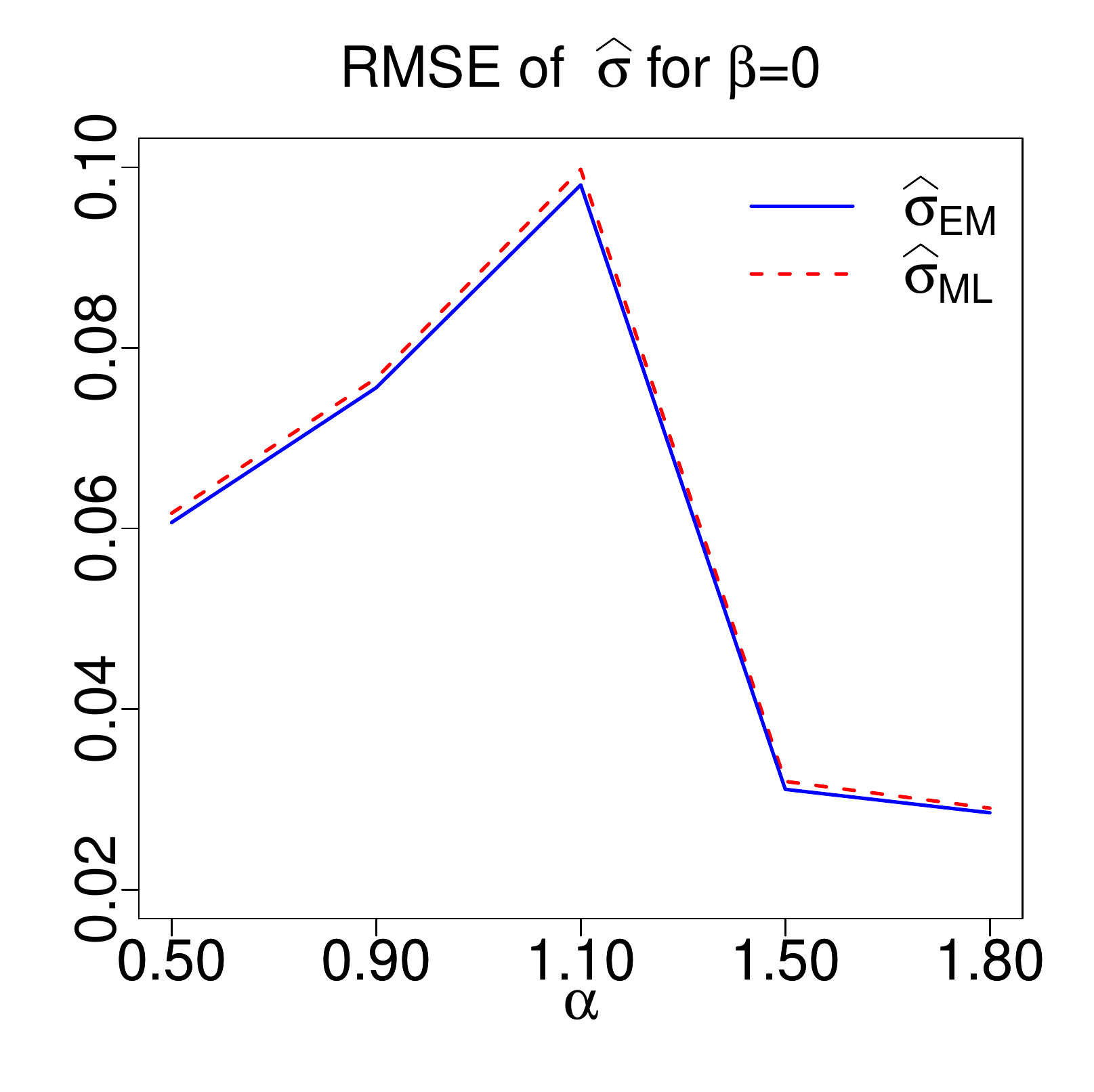}&
\includegraphics[width=40mm,height=40mm]{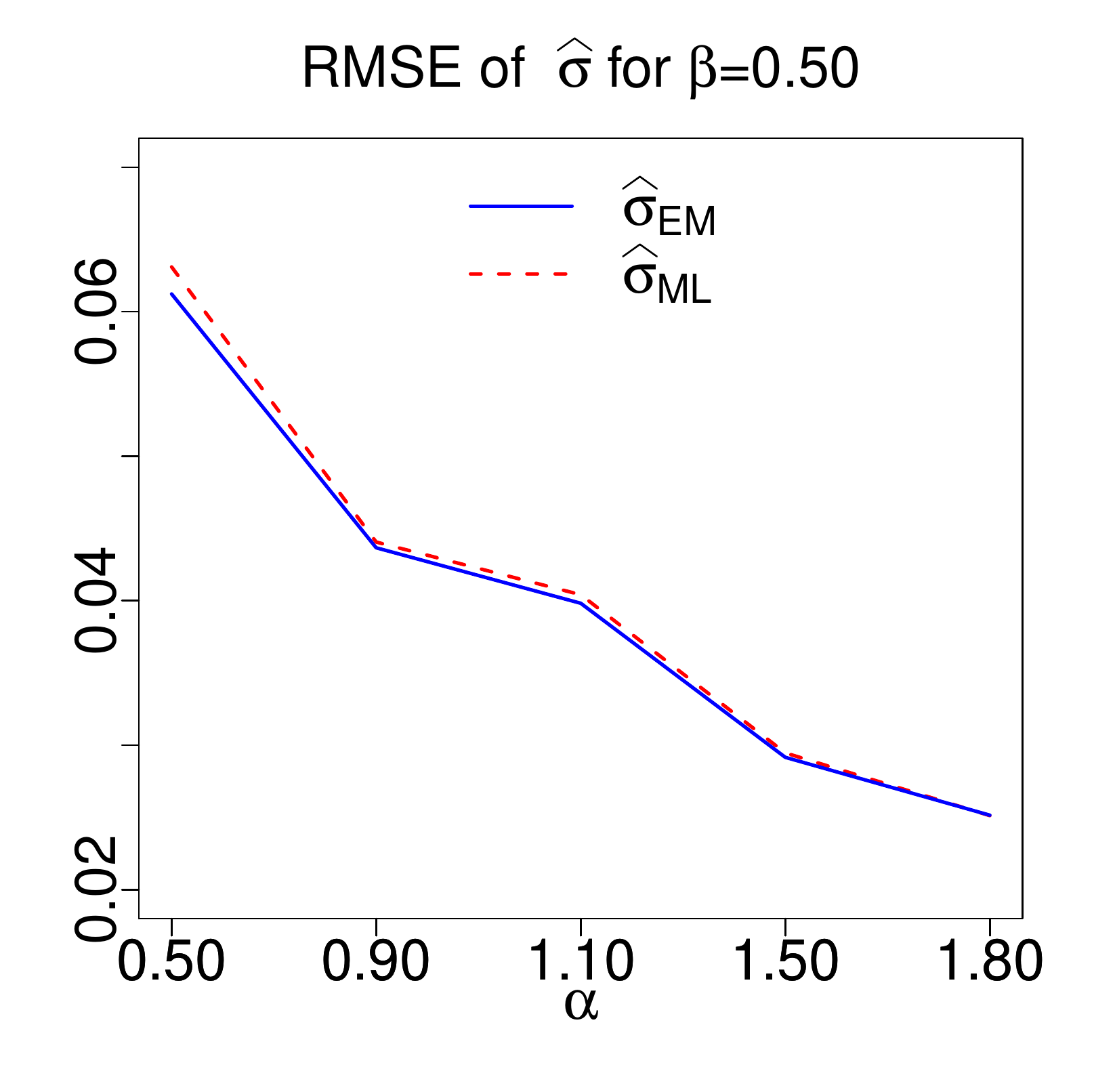}&
\includegraphics[width=40mm,height=40mm]{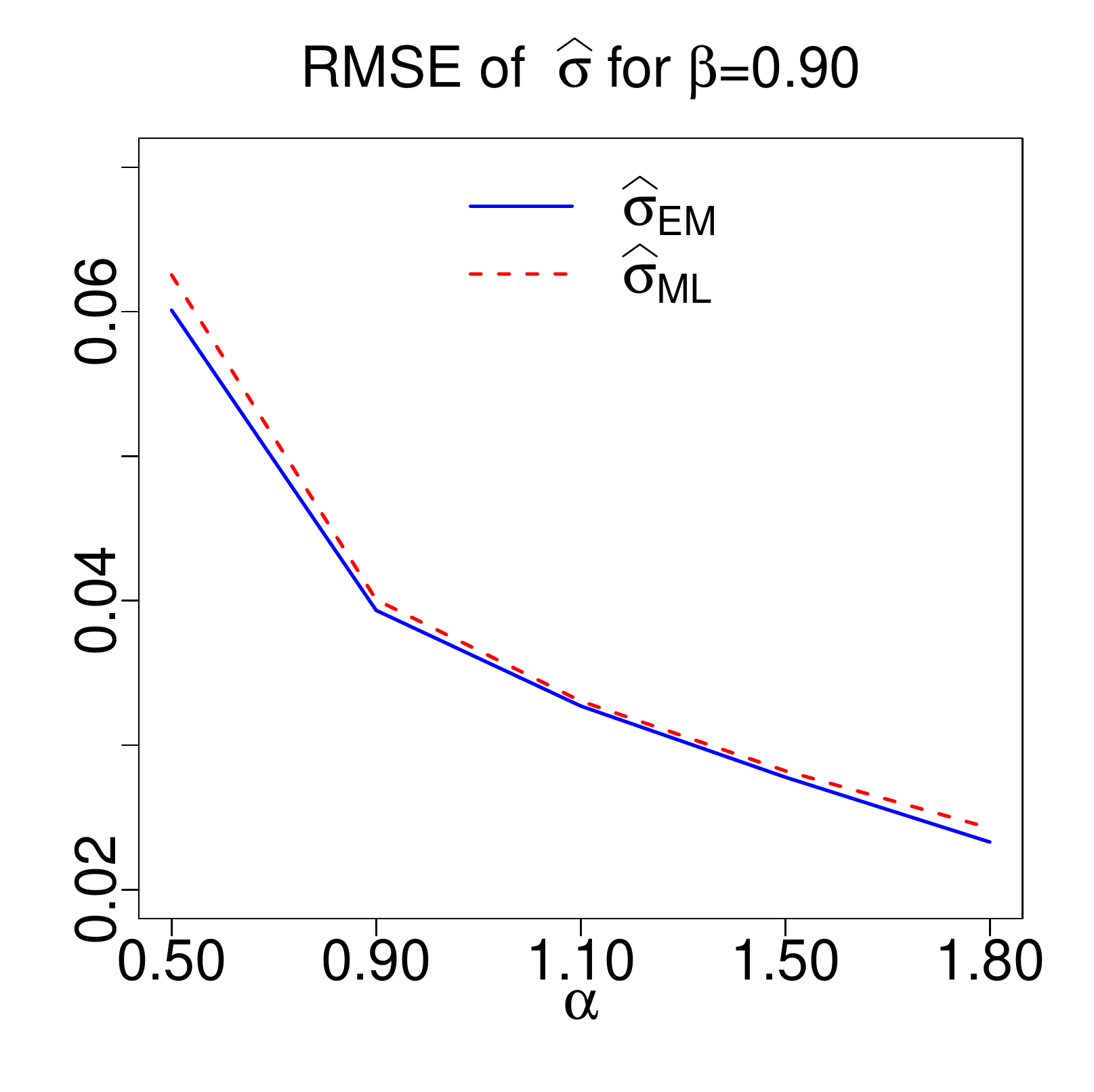}\\
\includegraphics[width=40mm,height=40mm]{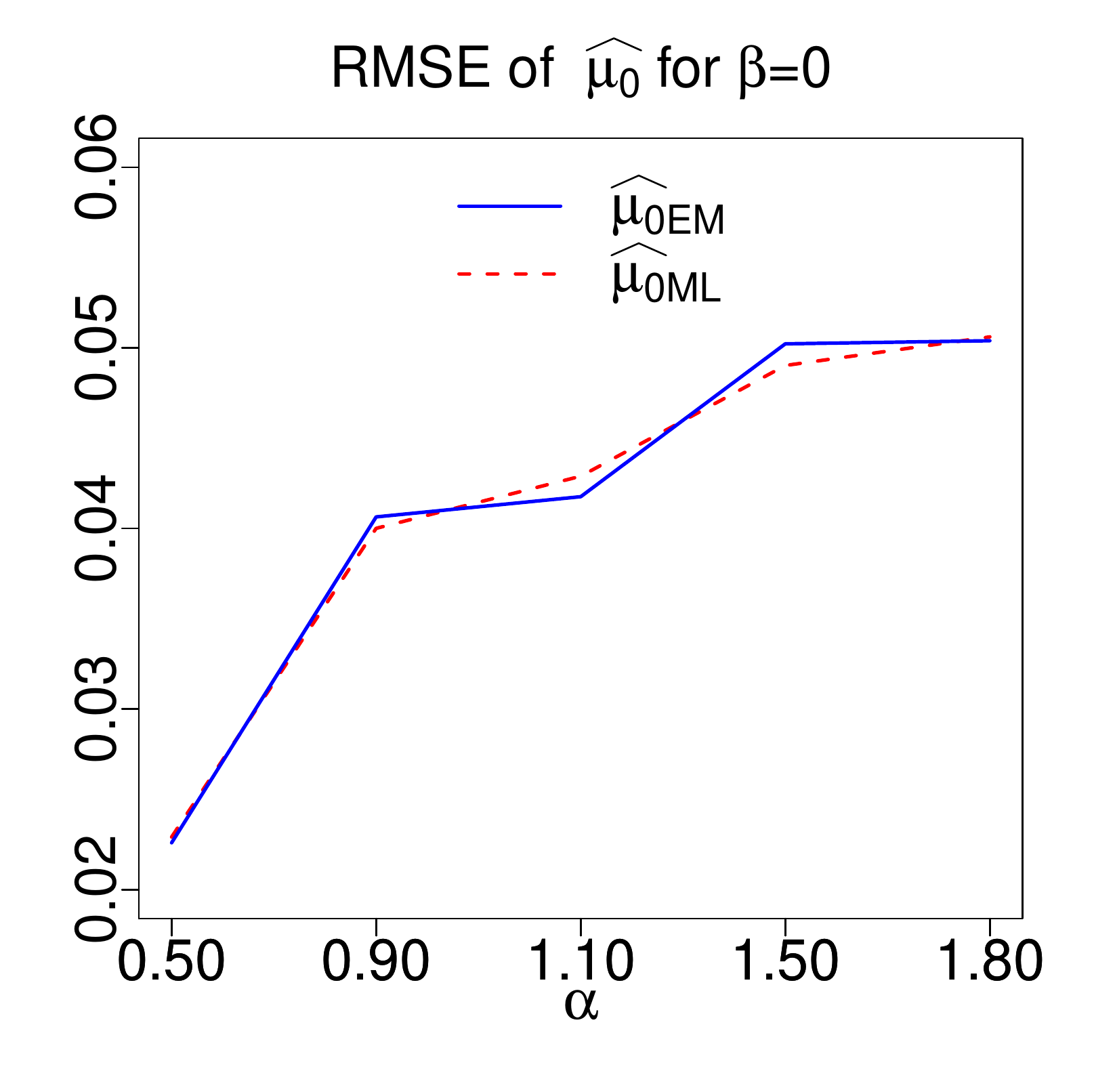}&
\includegraphics[width=40mm,height=40mm]{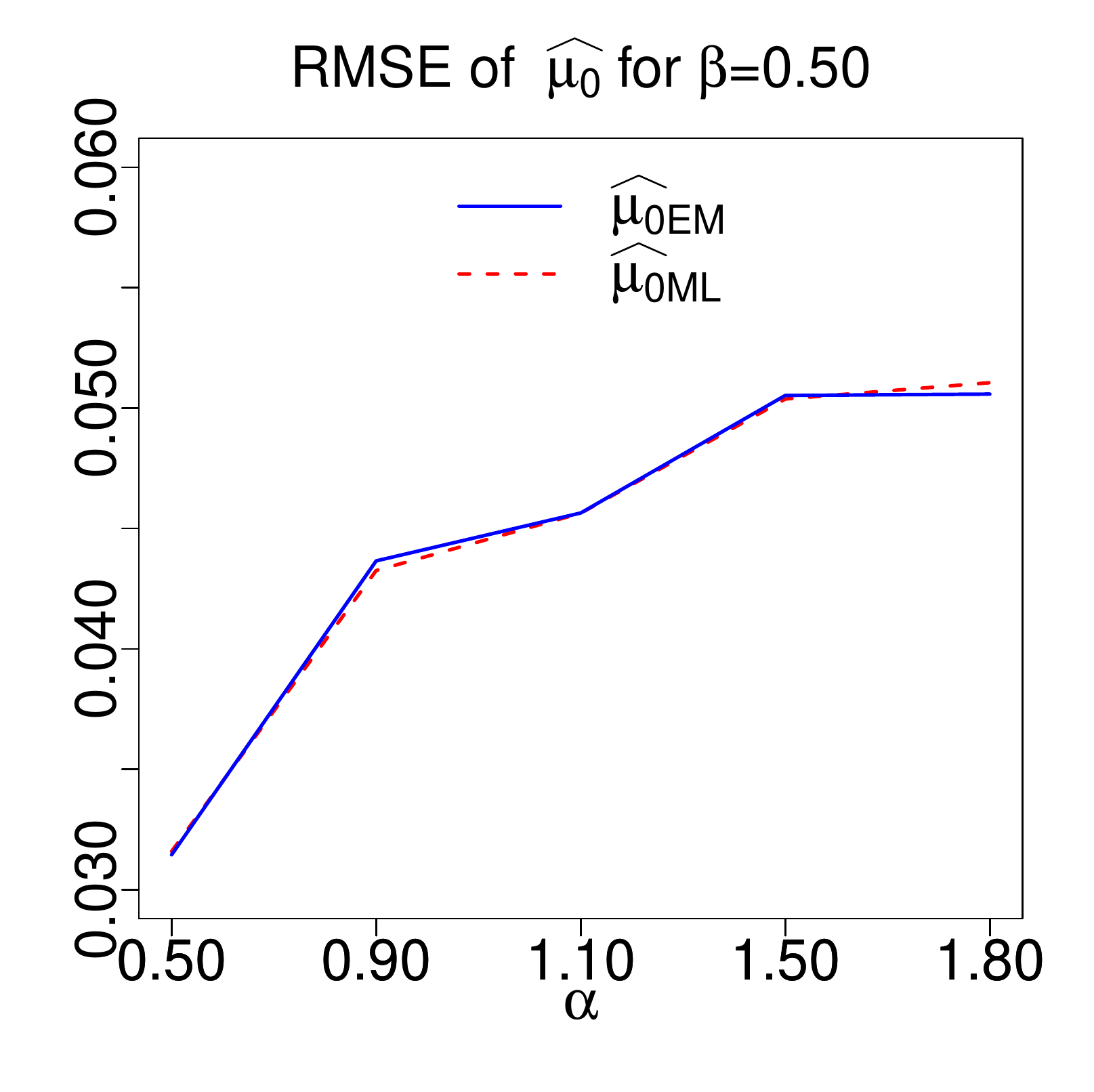}&
\includegraphics[width=40mm,height=40mm]{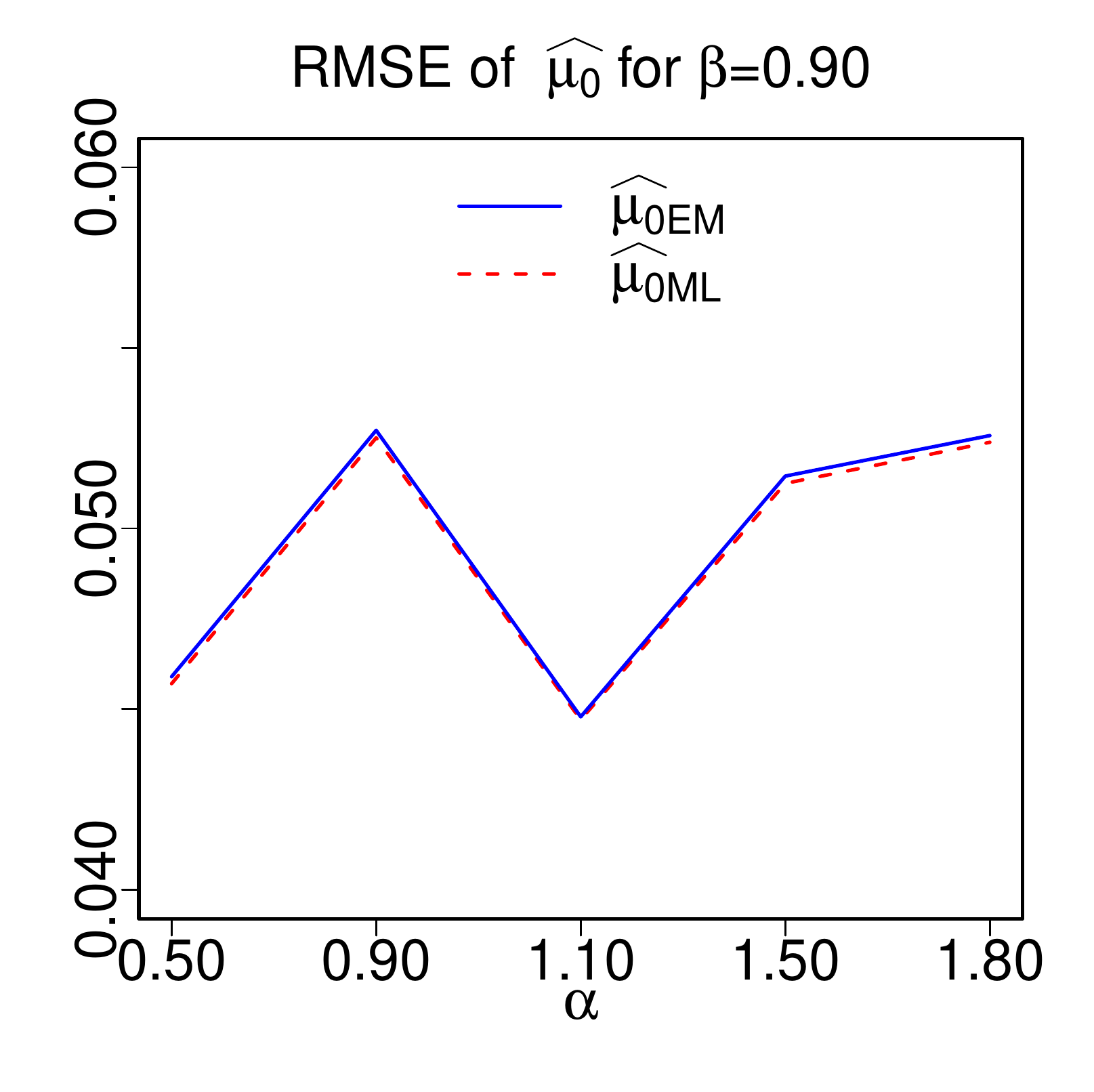}\\
\end{tabular}}
\caption{The RMSE of estimators obtained through the EM and ML approaches when $\sigma=0.5$ and $\mu_{0}=0$. In each sub-figure, the subscripts ML and EM indicate that the estimators $\hat{\alpha}$, $\hat{\beta}$, $\hat{\sigma}$, and $\widehat{\mu_{0}}$ are obtained using the EM algorithm (blue solid line) or the ML approach (red dashed line). Sub-figures in the first, second, and the third columns correspond to $\beta=0$, $\beta=0.50$, and $\beta=0.90$, respectively.}
\label{fig1}
\end{figure}

\begin{figure}
\resizebox{\textwidth}{!}
{\begin{tabular}{ccc}
\includegraphics[width=40mm,height=40mm]{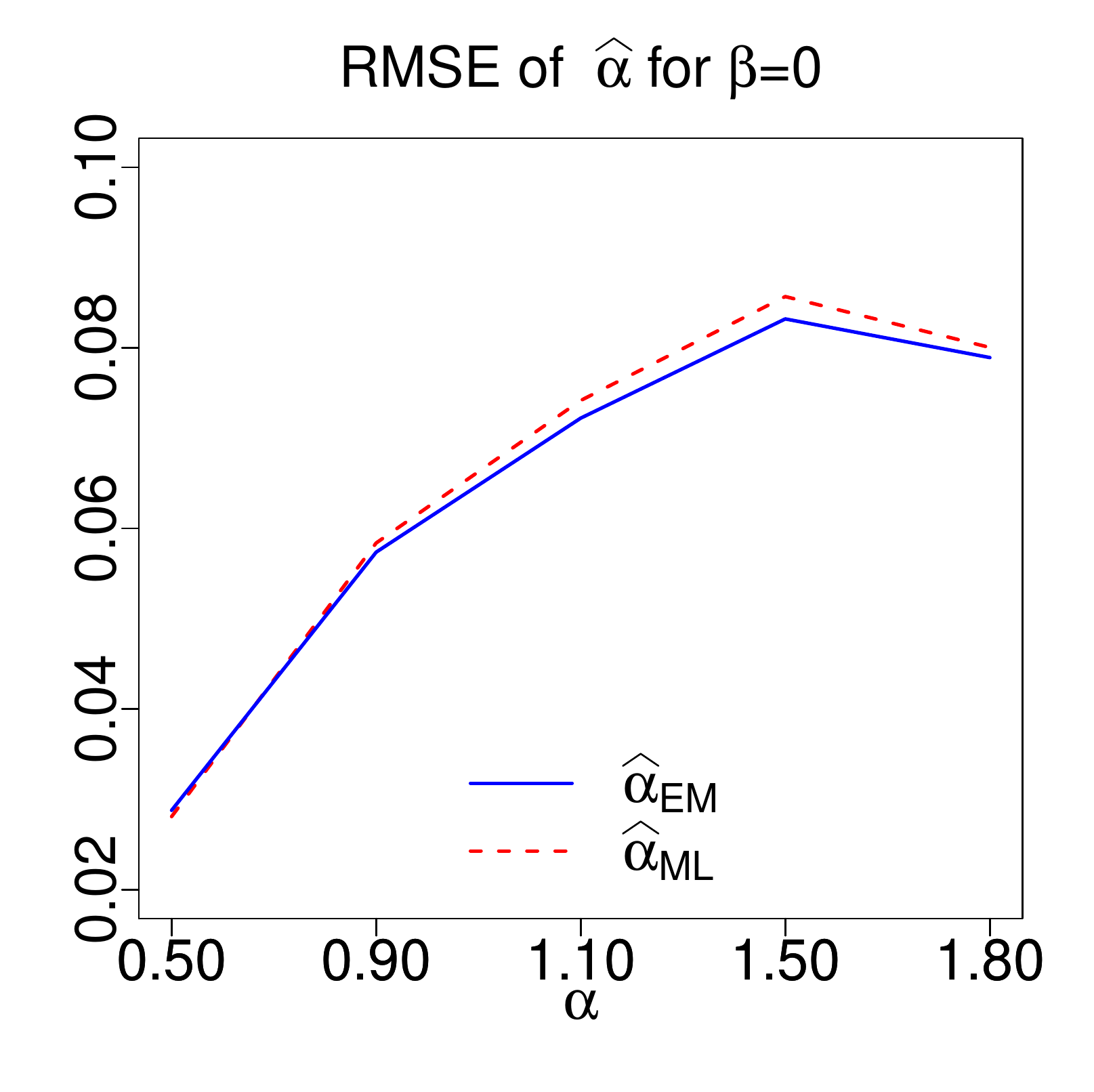}&
\includegraphics[width=40mm,height=40mm]{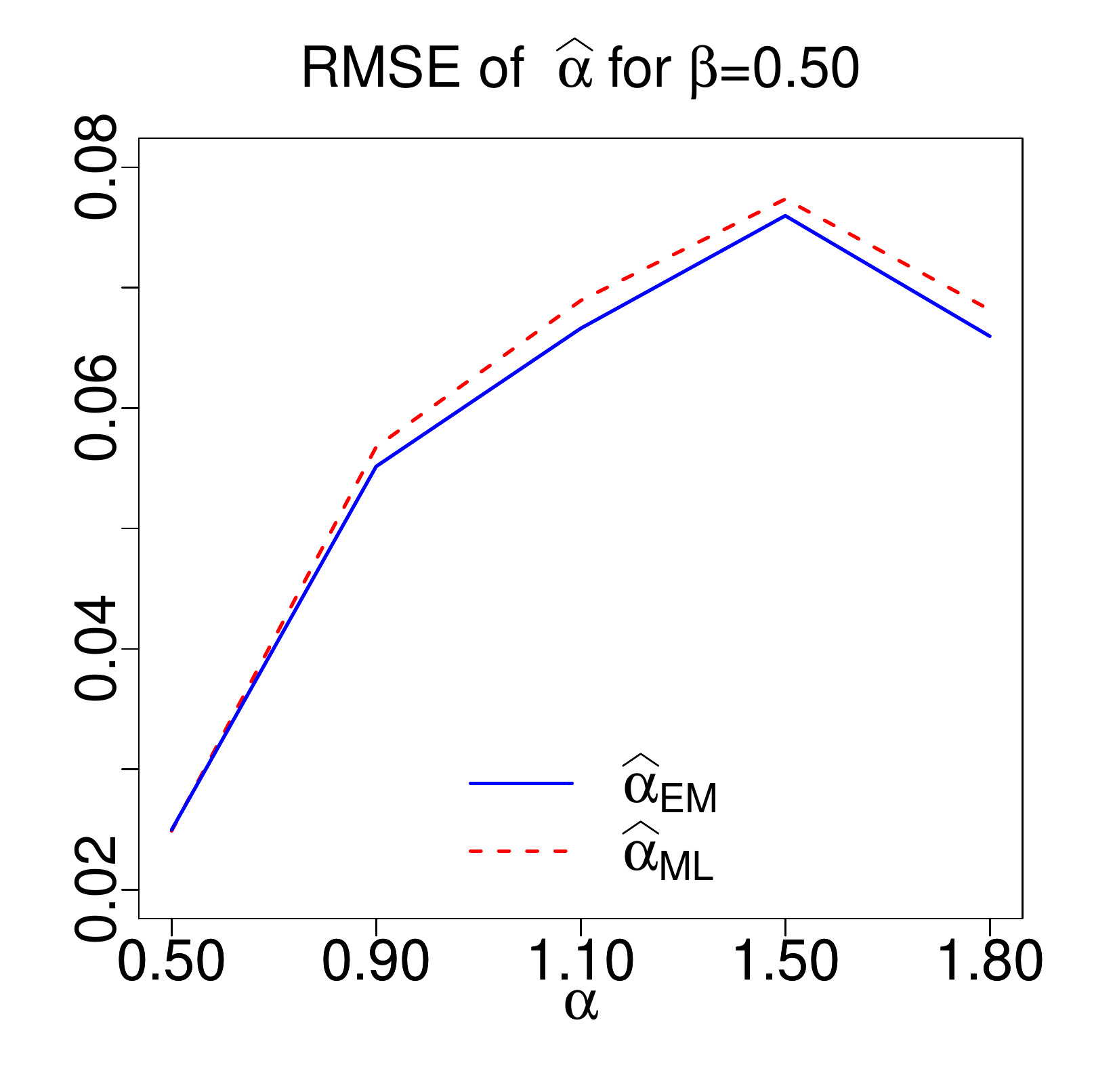}&
\includegraphics[width=40mm,height=40mm]{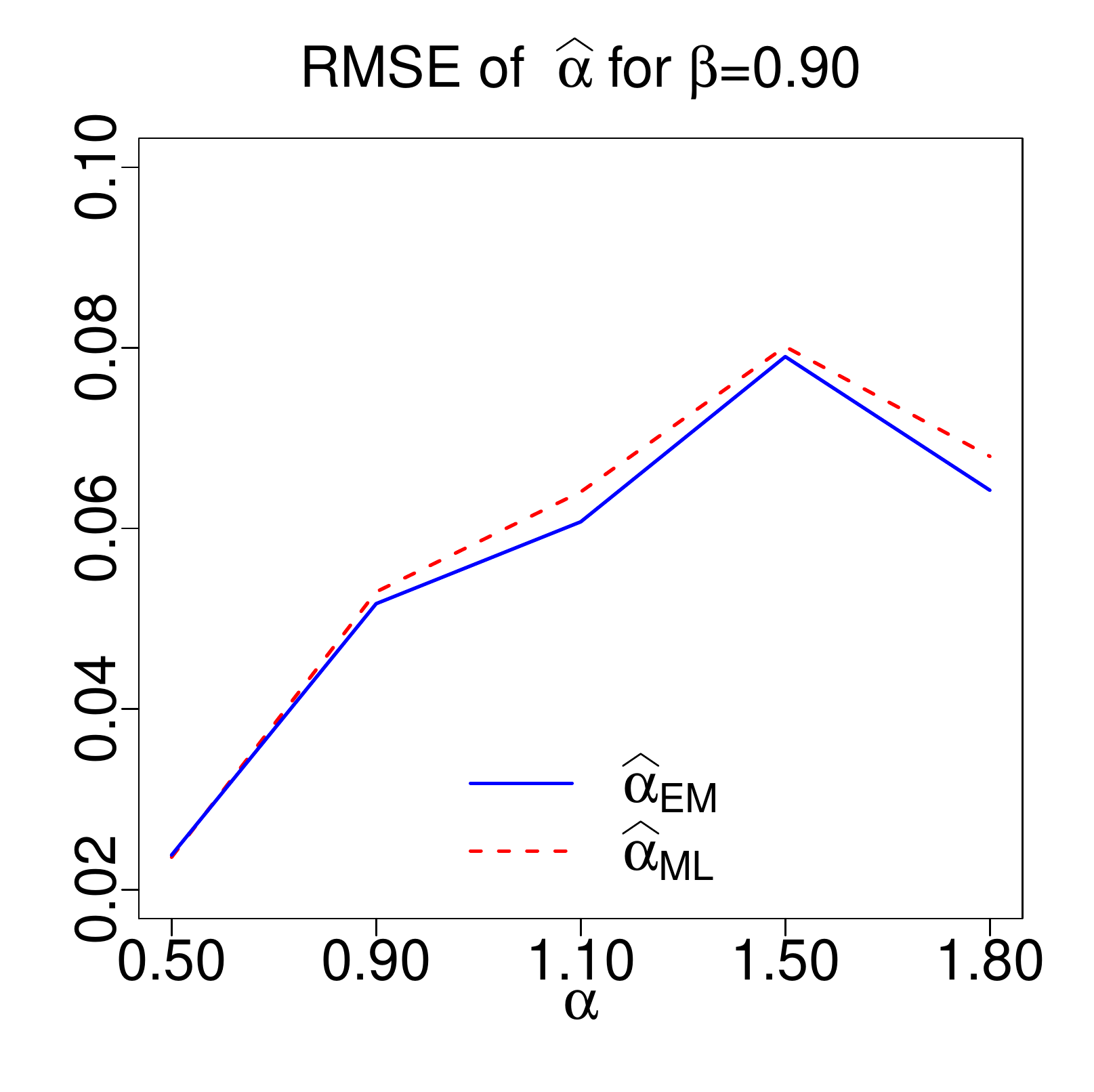}\\
\includegraphics[width=40mm,height=40mm]{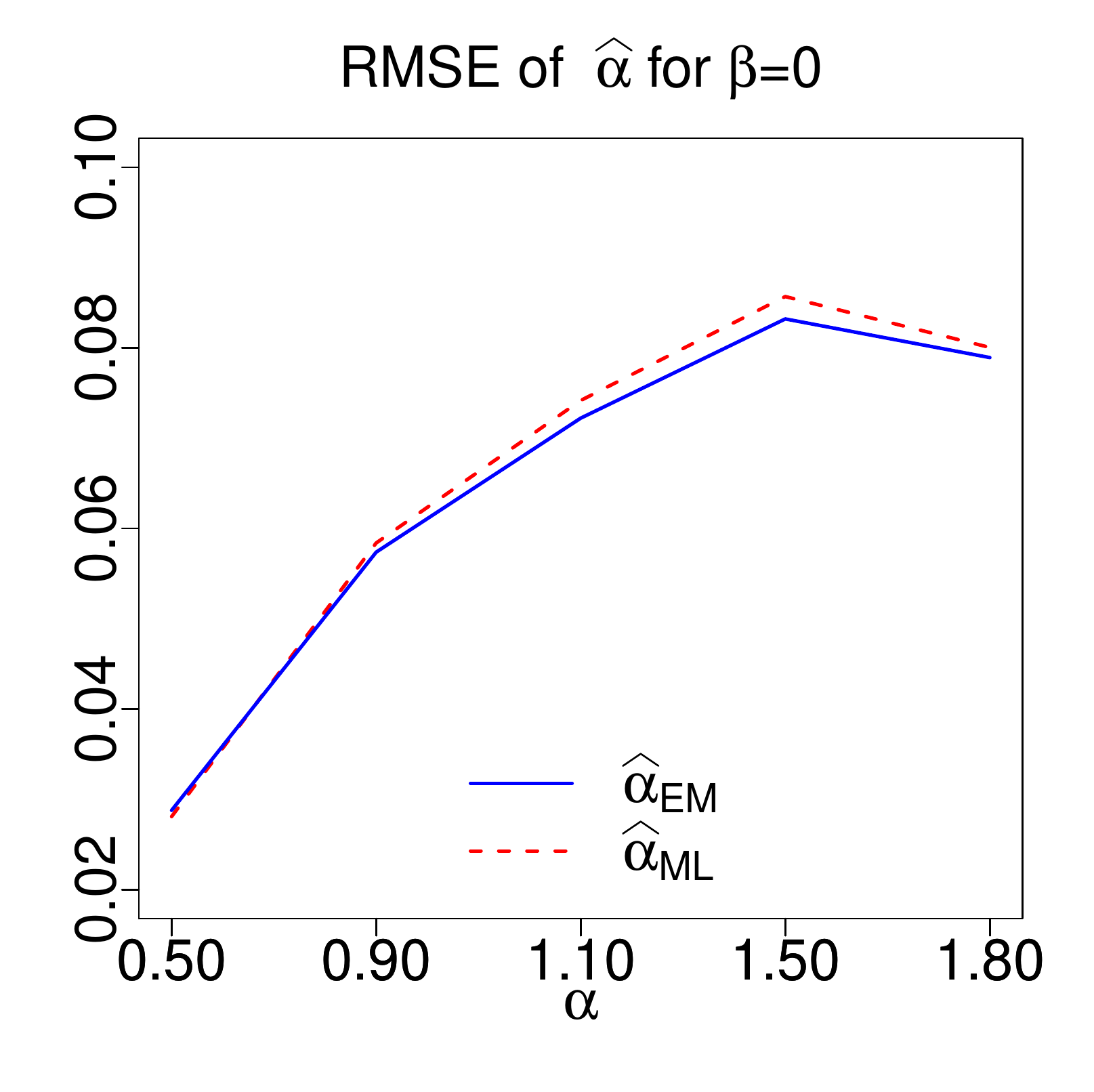}&
\includegraphics[width=40mm,height=40mm]{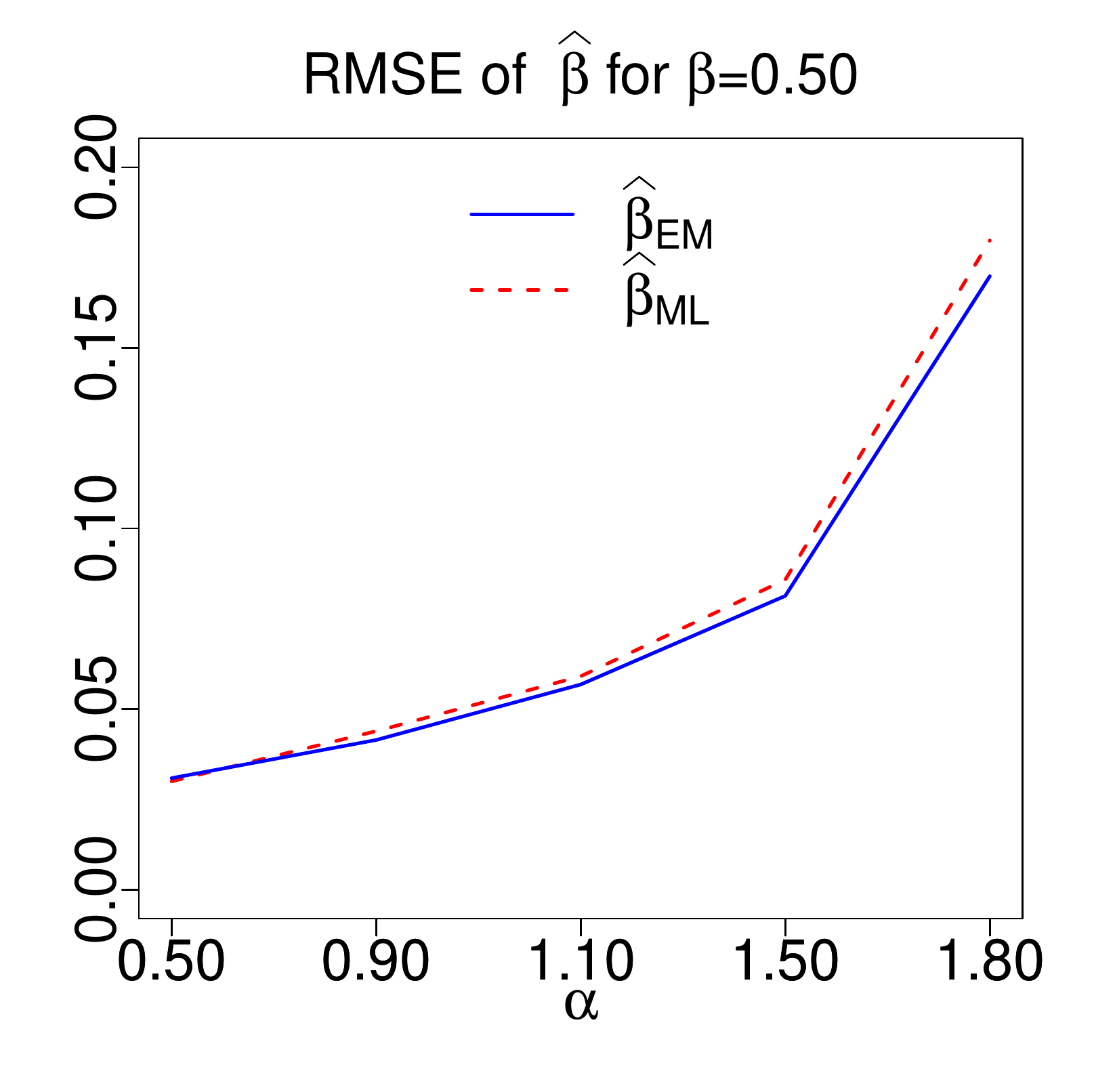}&
\includegraphics[width=40mm,height=40mm]{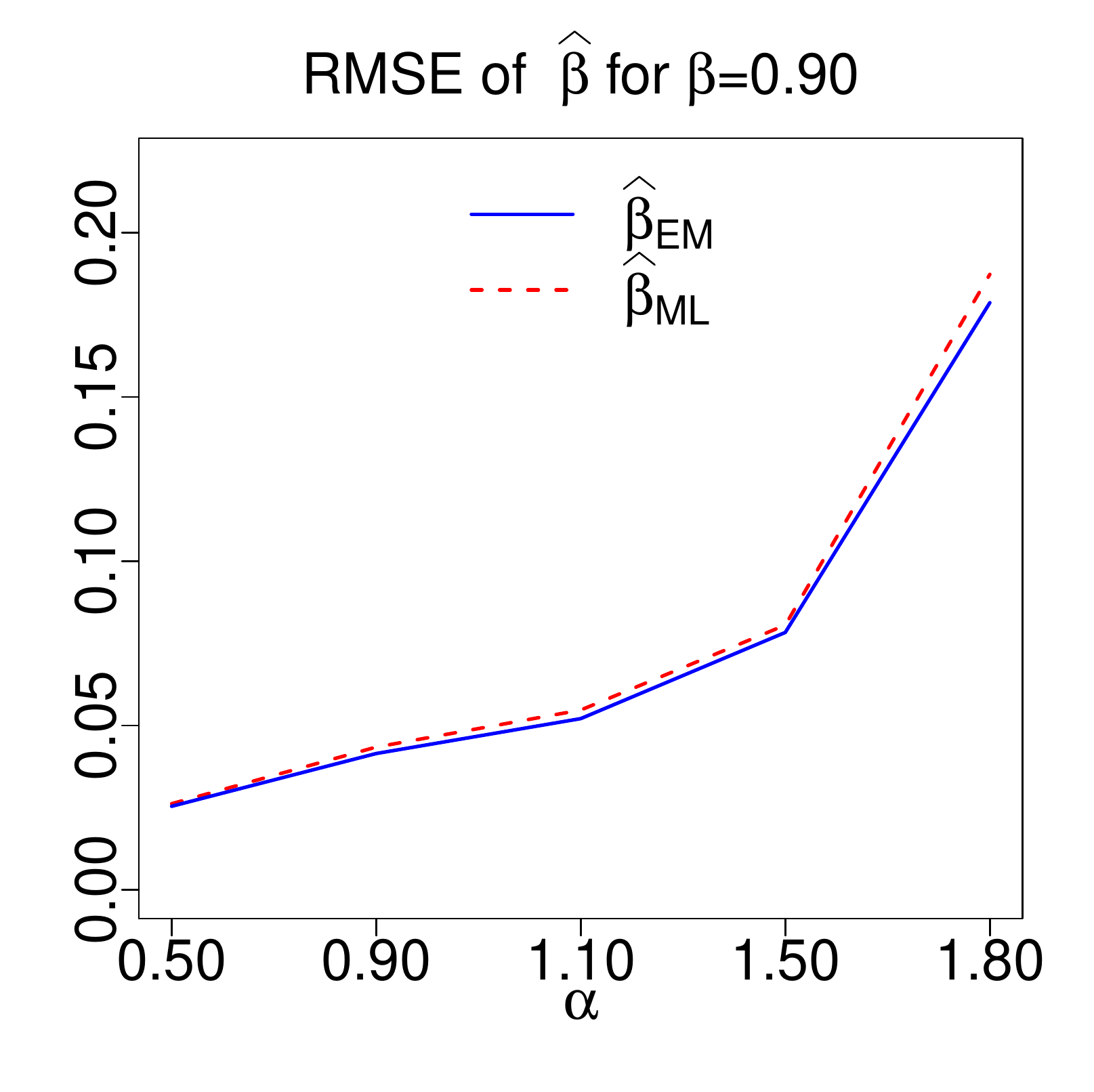}\\
\includegraphics[width=40mm,height=40mm]{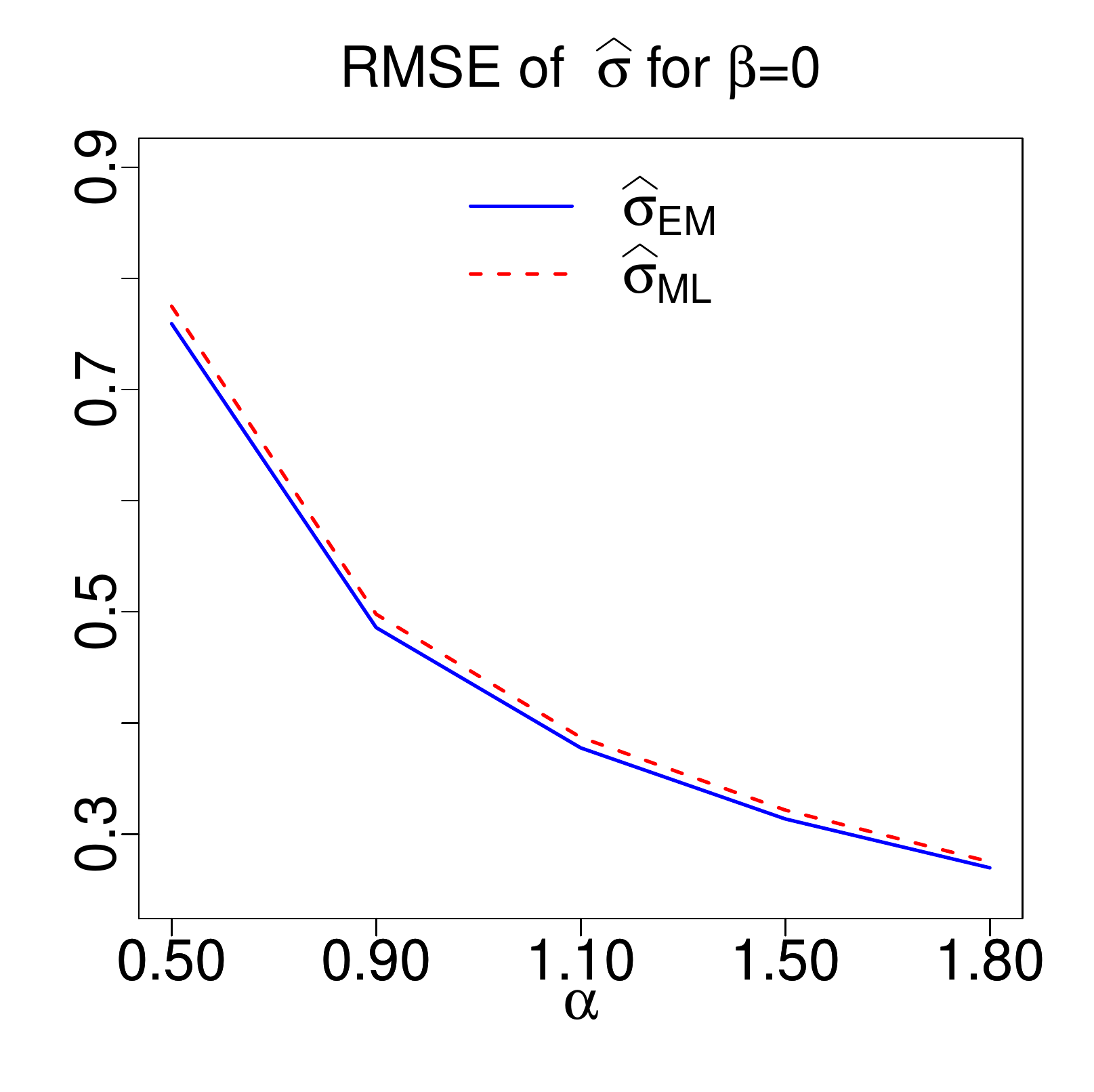}&
\includegraphics[width=40mm,height=40mm]{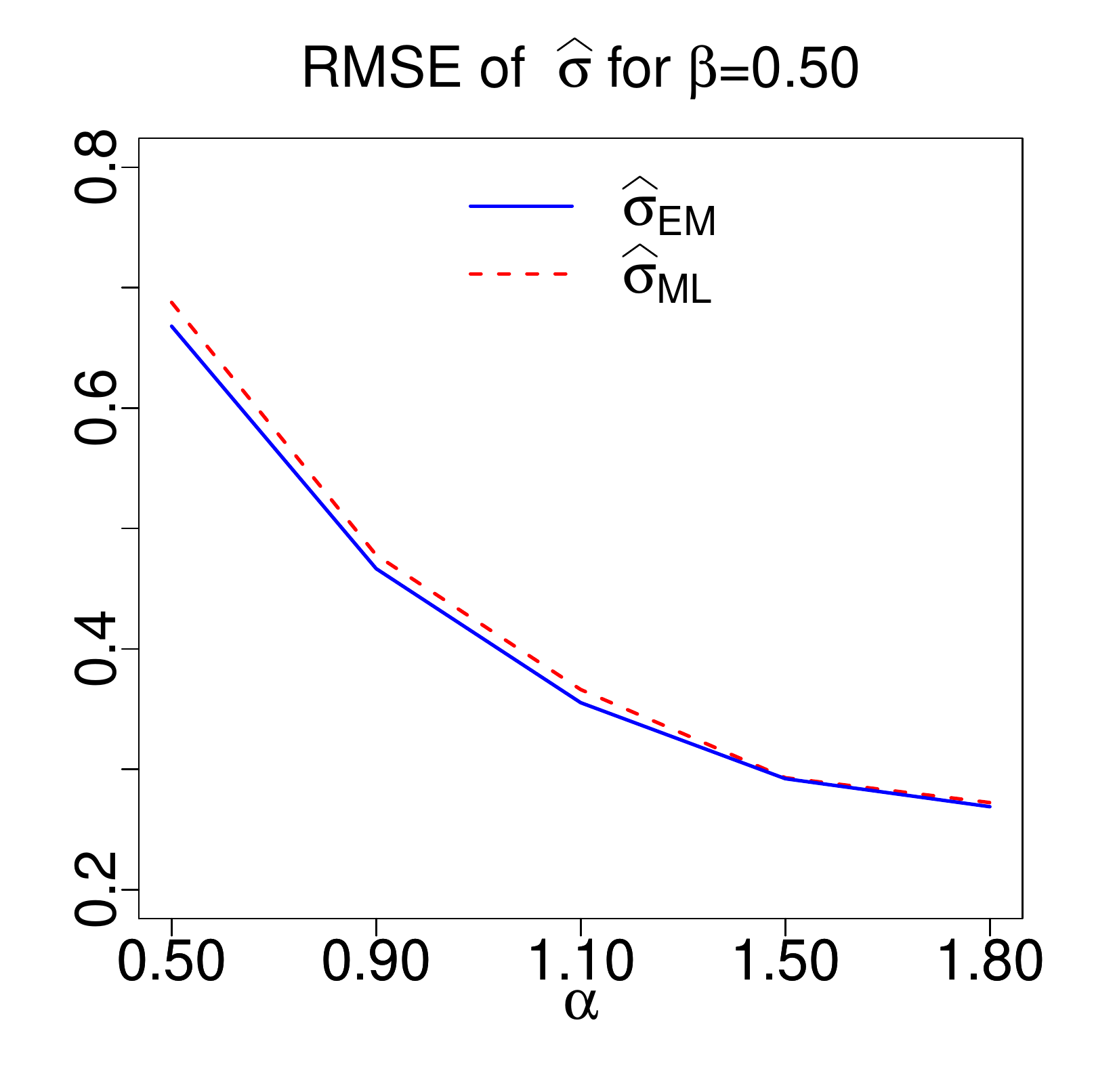}&
\includegraphics[width=40mm,height=40mm]{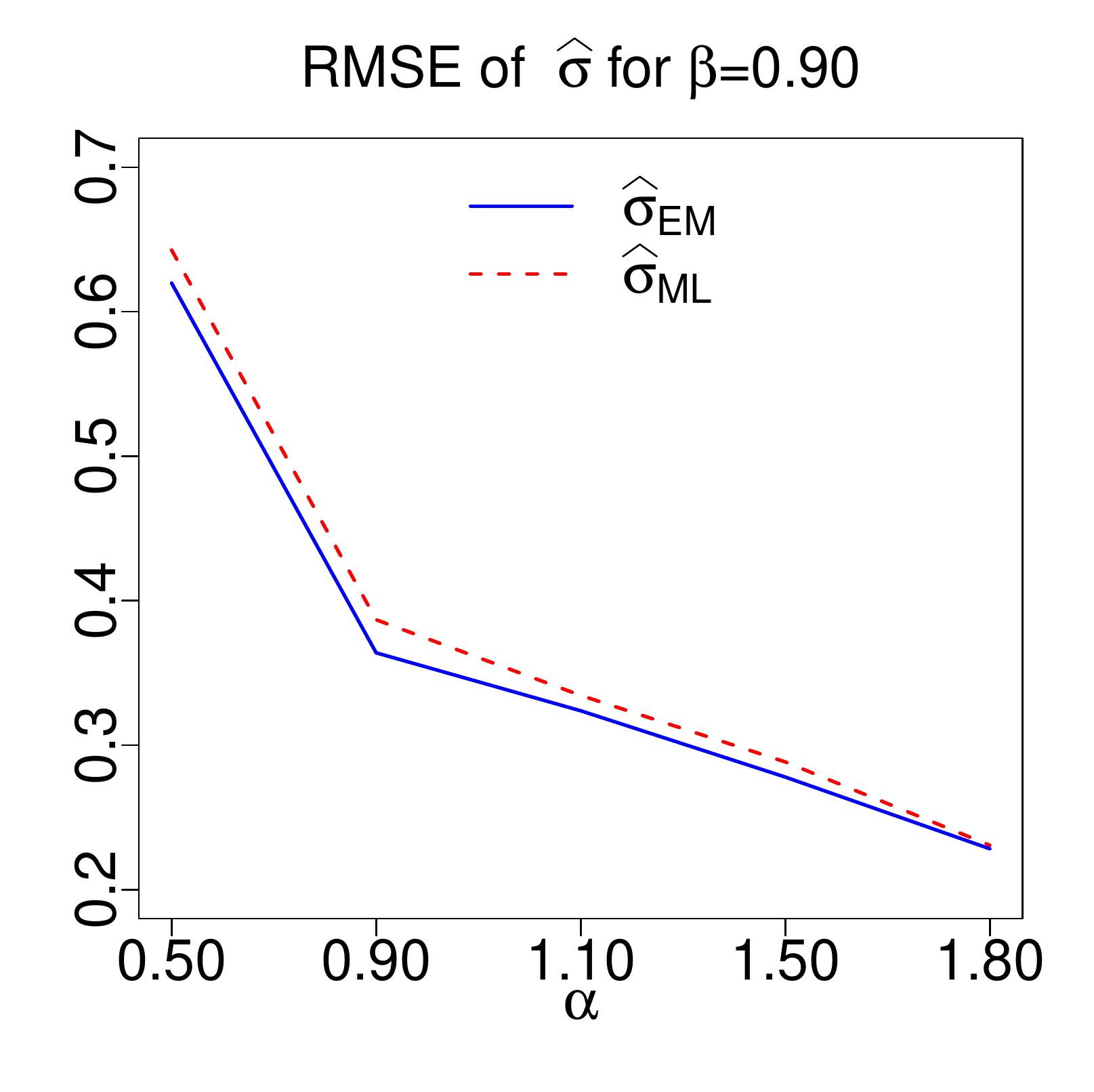}\\
\includegraphics[width=40mm,height=40mm]{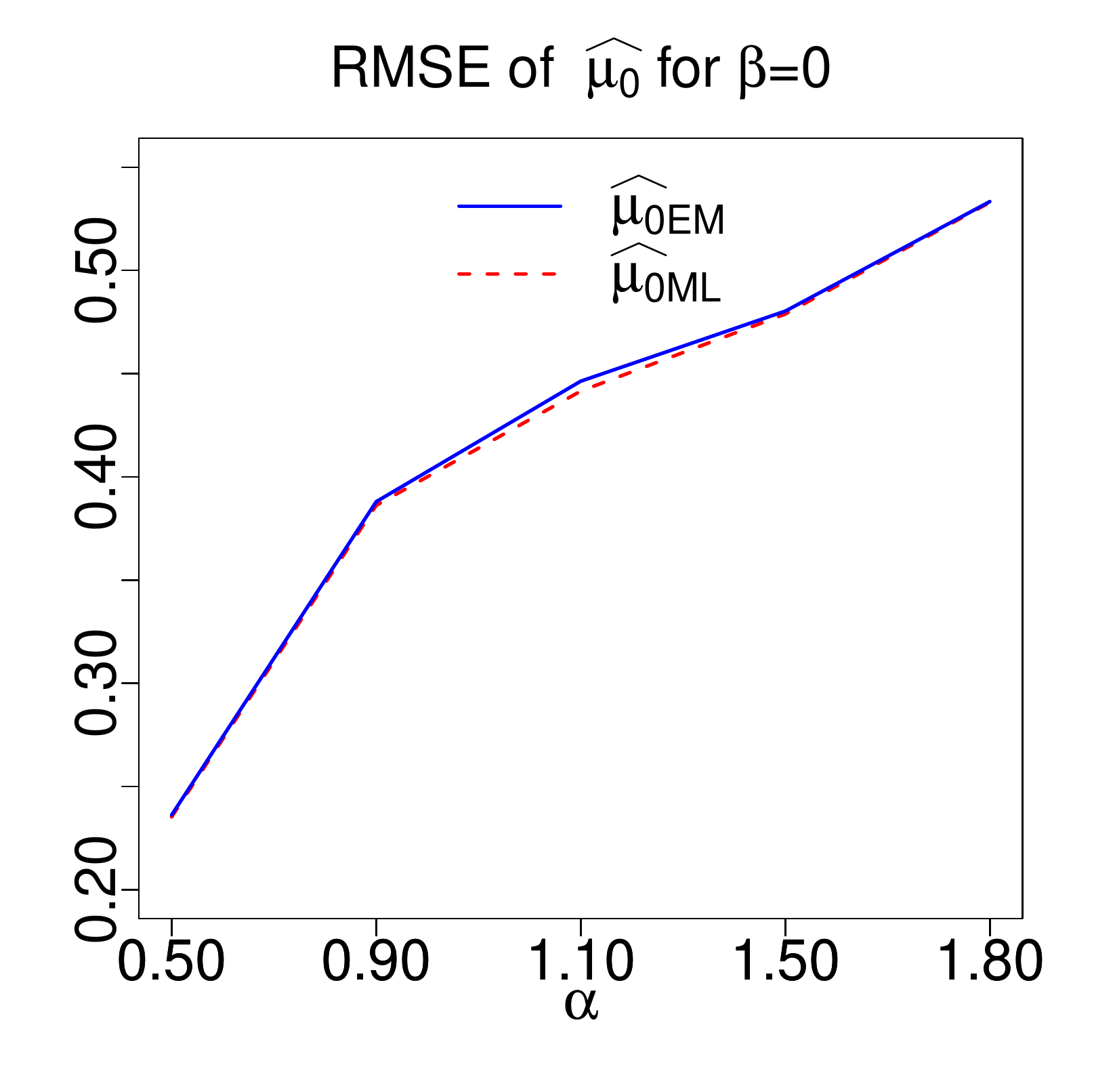}&
\includegraphics[width=40mm,height=40mm]{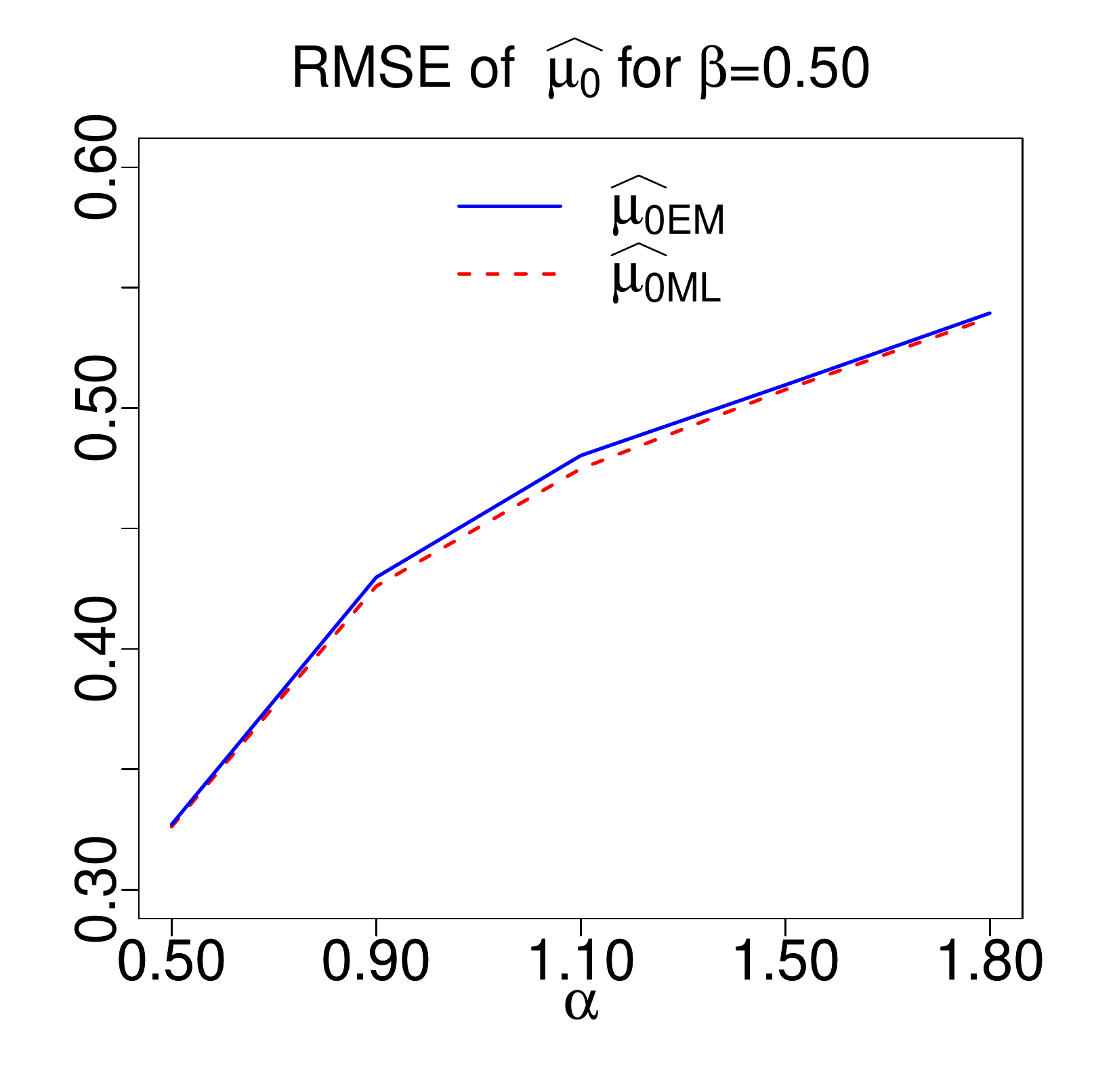}&
\includegraphics[width=40mm,height=40mm]{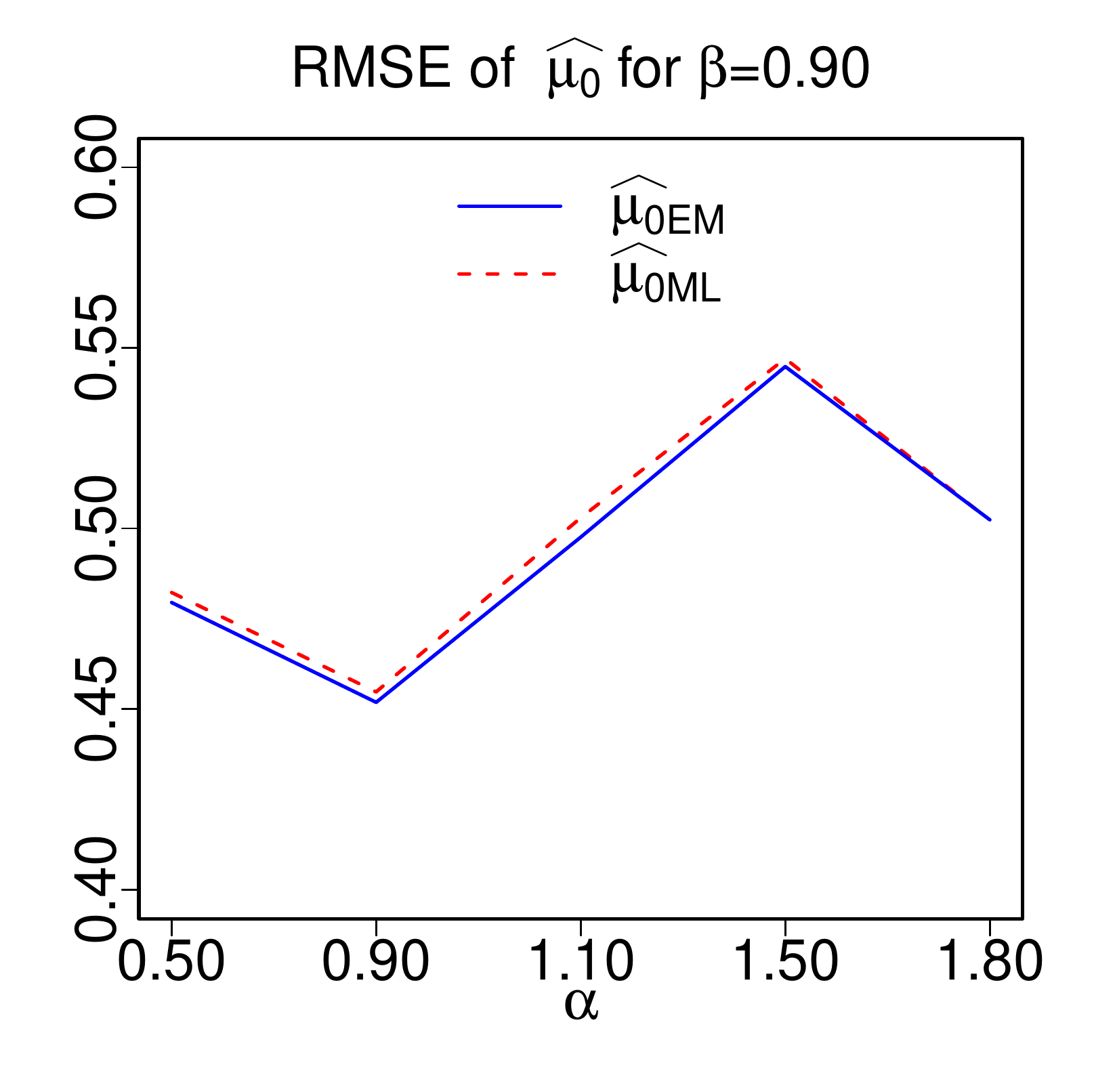}\\
\end{tabular}}
\caption{The RMSE of estimators obtained through the EM and ML approaches when $\sigma=5$ and $\mu_{0}=0$. In each sub-figure, the subscripts ML and EM indicate that the estimators $\hat{\alpha}$, $\hat{\beta}$, $\hat{\sigma}$, and $\widehat{\mu_{0}}$ are obtained using the EM algorithm (blue solid line) or the ML approach (red dashed line). Sub-figures in the first, second, and the third columns correspond to $\beta=0$, $\beta=0.50$, and $\beta=0.90$, respectively.}
\label{fig2}
\end{figure}

We apply the proposed EM and ML approaches to the 200 sets of samples of size 300. For the sake of simplicity, we restrict ourselves to the case of $\mu_{0}=0$. In each iteration, settings for the skewness and tail thickness parameters are: $\beta=0,0.5, 0.9$ and $\alpha=0.5, 0.9, 1.2, 1.5$. The results of simulations are shown in Figure \ref{fig1} and Figure \ref{fig2} for $\sigma=0.5$ and $\sigma=5$, respectively. During simulations, we set $M_0=20$, $M=40$, $N_0=100$, and $N=140$. This means that the proposed EM algorithm is run for 140 iterations of which the first 100 iterations are removed as burn-in period. Also, in each iteration, the CM-step is repeated for 40 times of which the average of the last 20 runs is considered as the updated tail thickness parameter.
\par Form Figure \ref{fig1} and Figure \ref{fig2}, in the sense of square root of the mean-squared error (RMSE) criterion, the following observations can be made.
\begin{itemize}
\item The EM-based estimator of the scale parameter outperforms the corresponding ML-based estimator when $\sigma=0.5$ (scale parameter is small).
\item The EM-based estimator of the tail thickness parameter, skewness, and scale parameters outperform the corresponding ML-based estimators when $\sigma=5$ (scale parameter is large). 
\end{itemize} 
It should be noted that the proposed EM algorithm shows better performance than the sample quantile (SQ, see \cite{McCulloch1986}) and empirical characteristic function (CF, see \cite{Kogon1998}) approaches, and so were eliminated by competitions. Also, as it is known, the EM algorithm cannot outperform the ML method. Sometimes, as noted above, the EM algorithm shows better performance than the ML method. This is because, $\mathsf{STABLE}$ computes the evaluated ML estimators not the exact ML ones. 
\subsection{Model validation via real data}
The $\alpha$-stable distribution is the most common used candidate for modelling the prices of
speculative assets such as stock returns, see \cite{buckle1995}, \cite{Mandelbrot}, \cite{Nolan2013}, and \cite{Rachev2005}. Here, we consider three sets of real data for illustrating an application of the $\alpha$-stable distribution. The data are daily price returns of the major European stock indices including Switzerland SMI, France CAC, and UK FTSE. Following common practice for daily closing prices, we consider the transformed prices for $n$ business days as  $(p_{t-1}-p_{t})/p_{t-1}$; for $t=2,\dots,n$. We then fitted the distribution to the transformed
data by the four methods including ML, EM, SQ, and CF. We obtained three sets of data from $\mathsf{datasets}$ package developed for $\mathsf{R}$ (R Core Team \cite{rteam}) environment which include $n=1860$ observations. The results after applying above four approaches are given in Table \ref{tab1}. It should be noted that, for implementing the EM algorithm, the started $\alpha$, $\beta$, $\sigma$, and $\mu_0$ are well away from the estimated values. 

\begin{figure}
\resizebox{\textwidth}{!}
{\begin{tabular}{cc}
\includegraphics[width=40mm,height=40mm]{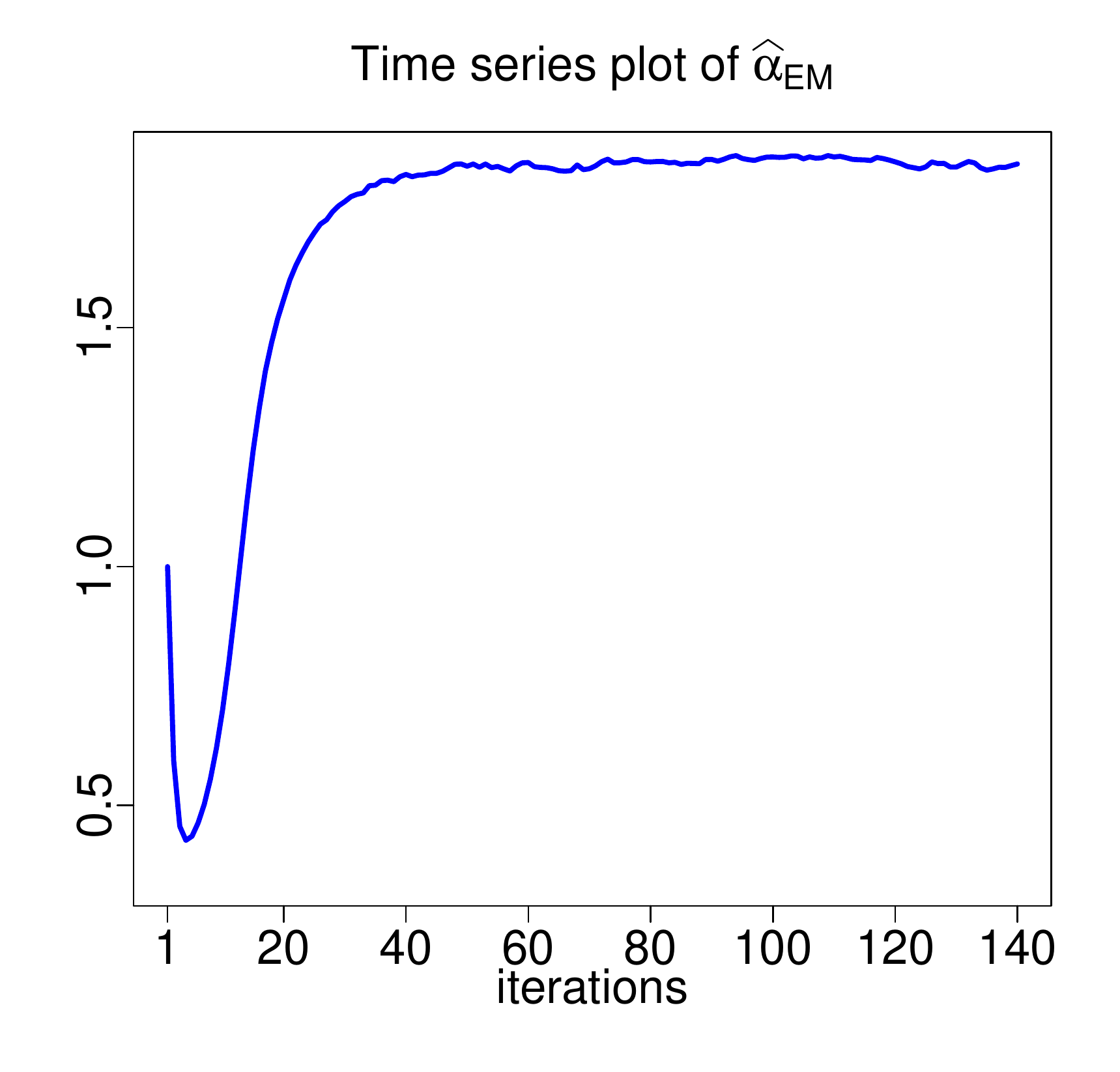}&
\includegraphics[width=40mm,height=40mm]{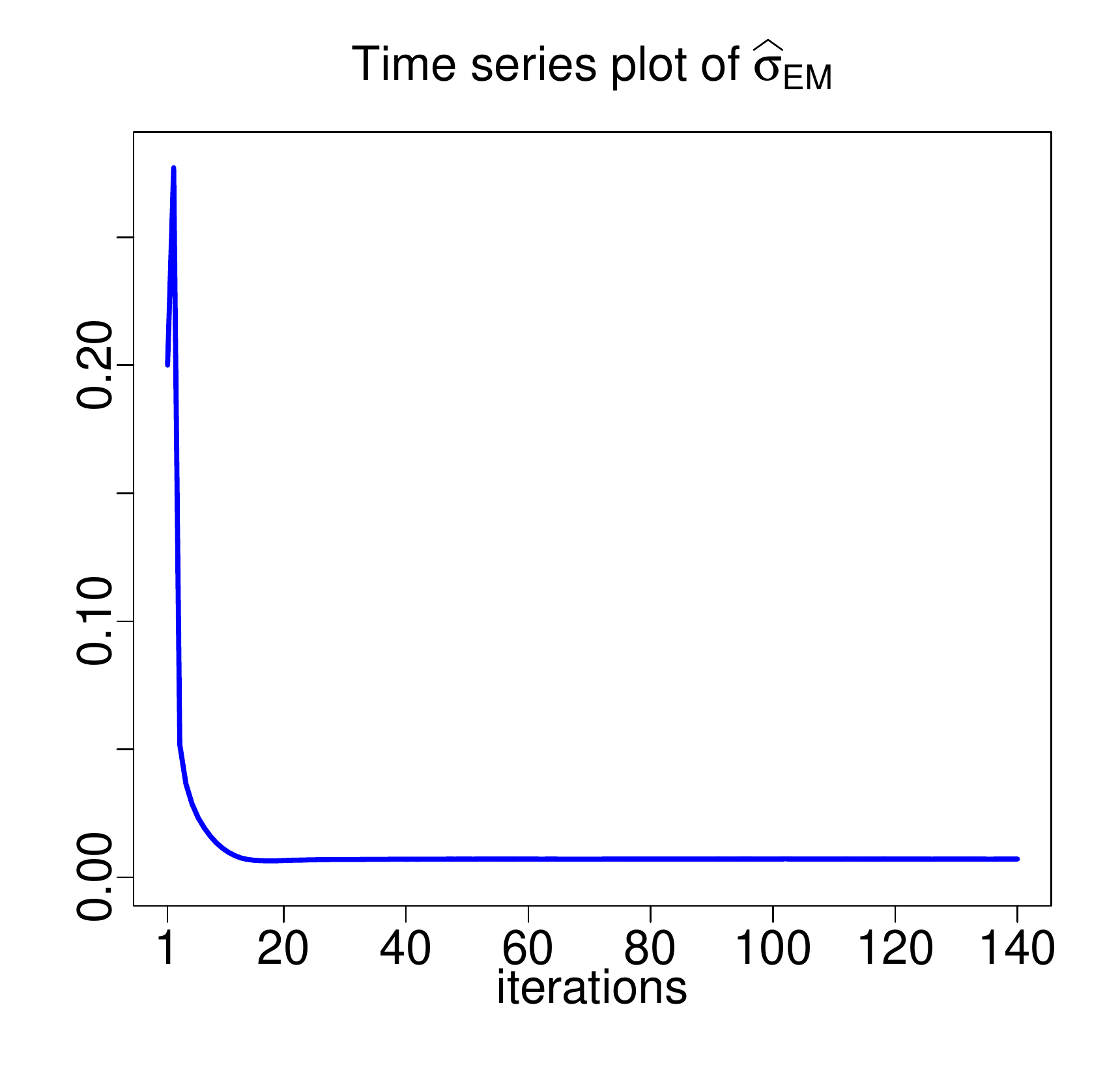}\\
\includegraphics[width=40mm,height=40mm]{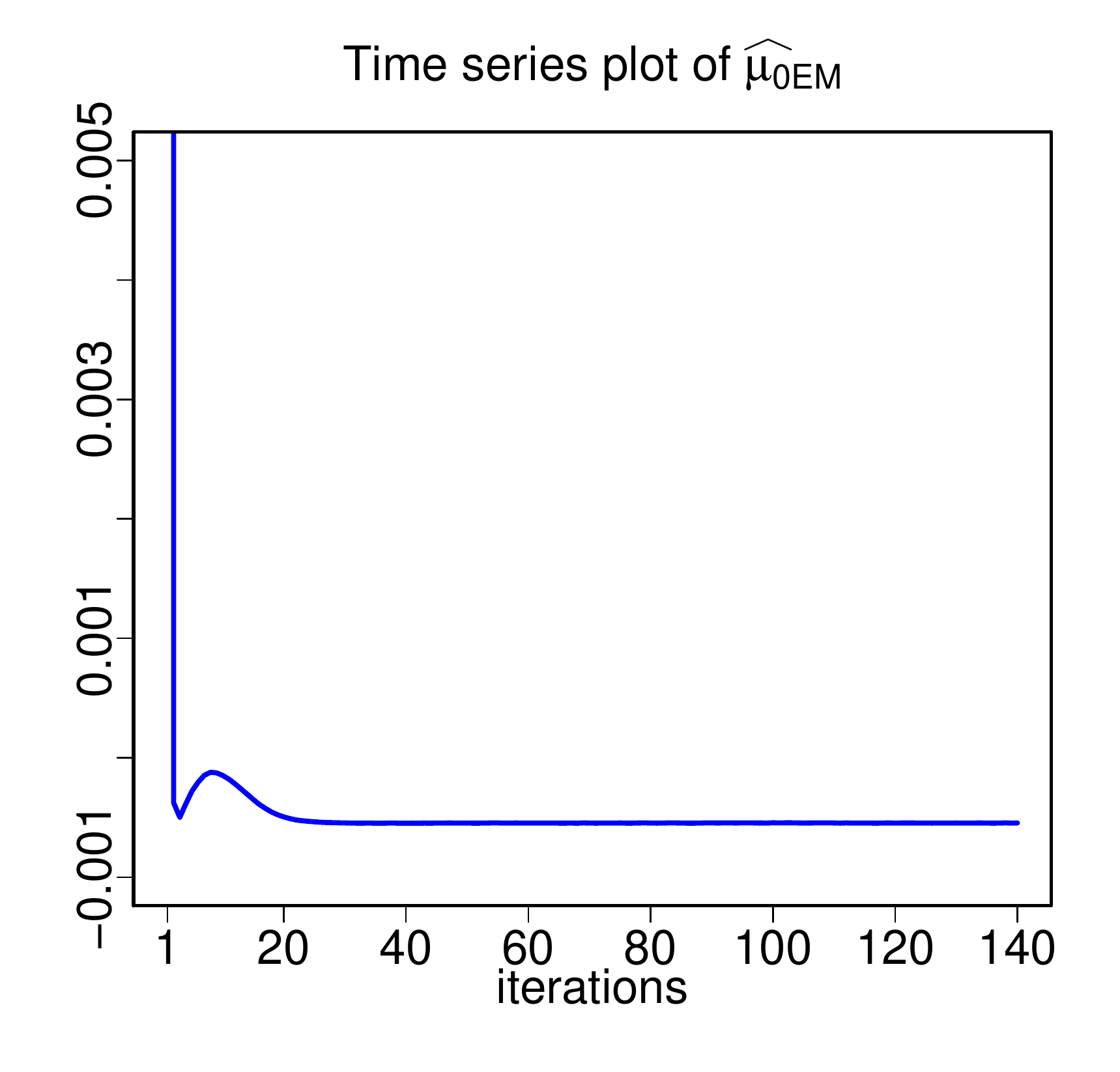}&
\includegraphics[width=40mm,height=40mm]{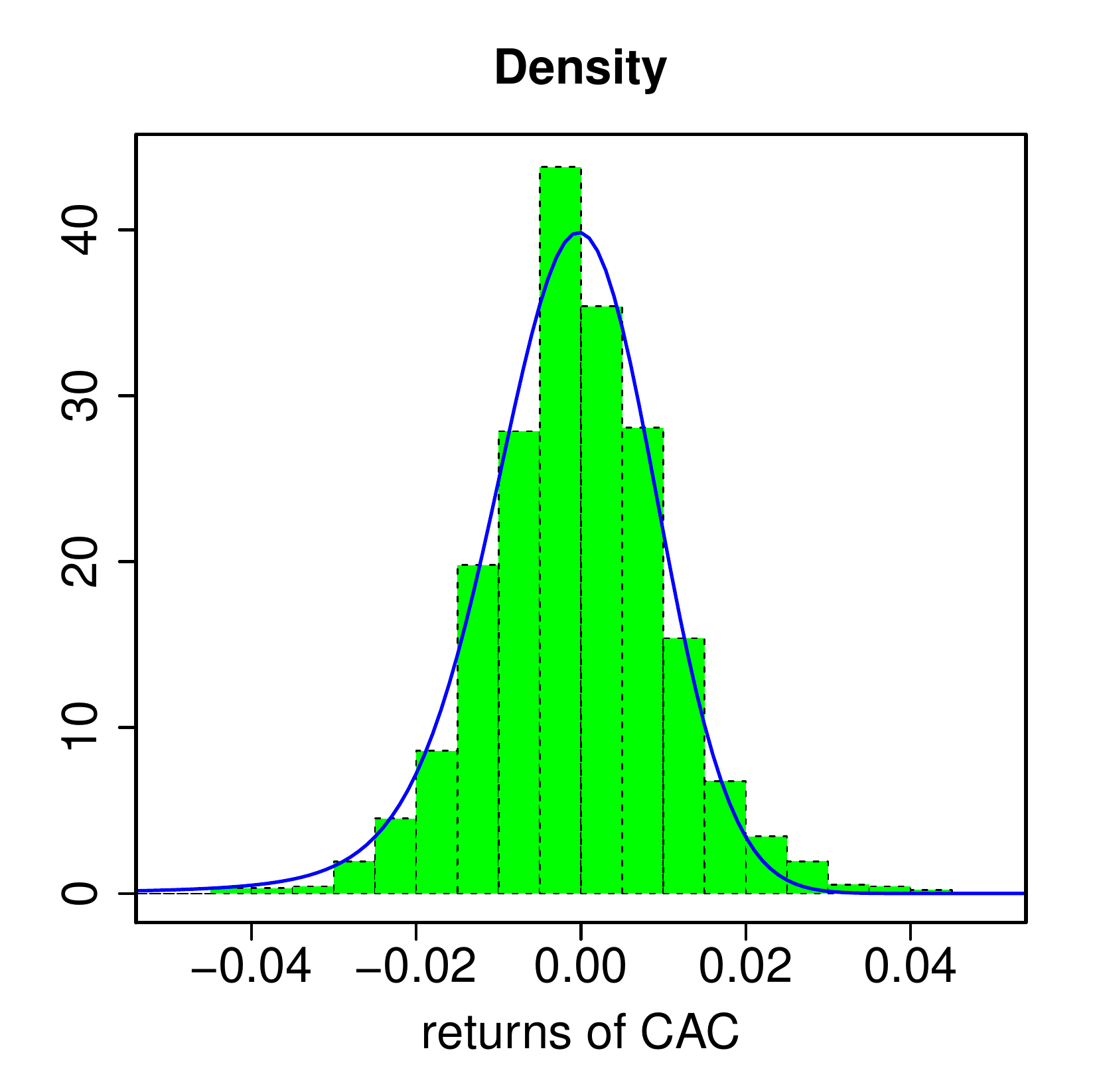}\\
\end{tabular}}
\caption{Time series plots of the updated parameters and histogram of 1859 CAC returns. To implement the proposed EM algorithm the initial values for $\alpha$, $\beta$, $\sigma$, and $\mu_{0}$ are 0.8, 0, 0.25, and 0.25, respectively. Fitted $S_{0}\bigl(\hat{\alpha}_{EM},\hat{\beta}_{EM},\hat{\sigma}_{EM},\widehat{\mu_{0}}_{EM})$ distribution captures clearly the histogram of data.}
\label{fig3}
\end{figure}

\begin{table}[h!]
\center
\caption{Estimated parameters using the proposed EM and ML methods for five sets of real data. The log-likelihood and Kolmogorov-Smirnov (KS)
statistics are given.} 
\begin{tabular}{cccccccc} 
\cline{1-8} 
& & \multicolumn{4}{c}{Estimated parameter}\\ \cline{3-6} 
Data set&Approach&$\hat{\alpha}$&$\hat{\beta}$&$\hat{\sigma}$&$\hat{\mu_{0}}$&Log-likelihood & KS\\ \cline{1-8} 
\cline{1-8}
\multicolumn{1}{c}{\multirow{2}{*}{SMI}}& 
\multicolumn{1}{c}{{{EM}}} &1.76460&0.15189&0.00541&-0.00106&{{6168.845}}&0.03335\\ 
&\multicolumn{1}{c}{{{ML}}}&1.74707&0.19603&0.00543&-0.00126&6168.528&0.02531\\ 
&\multicolumn{1}{c}{{{SQ}}}&1.60081&0.06607&0.00512&-0.00096&6161.445&0.03350\\ 
&\multicolumn{1}{c}{{{CF}}}& 1.81778&0.23881&0.00549&-0.00128&6167.198&0.02456\\ 
\cline{1-8}
\\
\cline{1-8}
\multicolumn{1}{c}{\multirow{2}{*}{CAC}}& 
\multicolumn{1}{c}{{{EM}}} &1.84712&0.04423&0.00707&-0.00054&5780.248&{{0.03018}}\\
&\multicolumn{1}{c}{{{ML}}}&1.86714&0.08863&0.00713&-0.00062&5780.415&0.03256 \\ 
&\multicolumn{1}{c}{{{SQ}}}&1.76230&-0.09998&0.00685&0.00009&5773.361&0.04129\\ 
&\multicolumn{1}{c}{{{CF}}}&1.90333&-0.08305&0.00712&-0.00050&5779.068& 0.03151\\ 
\cline{1-8}
\\
\cline{1-8}
\multicolumn{1}{c}{\multirow{2}{*}{FTSE}}& 
\multicolumn{1}{c}{{{EM}}} &1.86871&0.03927&0.00514&-0.00045&6396.488&{{0.02097}}\\ 
&\multicolumn{1}{c}{{{ML}}}&1.86640&0.07575&0.00510&-0.00049&6396.572&0.02237\\ 
&\multicolumn{1}{c}{{{SQ}}}&1.76710& -0.05947&0.00498&-0.00003&6389.608&0.03927\\ 
&\multicolumn{1}{c}{{{CF}}}&1.90178&  0.08753&  0.00512& -0.00049&6395.707&0.02231\\ 
\cline{1-8}
\\
\end{tabular} 
\label{tab1}
\end{table} 
It is clear from Table \ref{tab1} that the proposed EM algorithm outperforms the M, SQ, and CF approaches with respect to either the log-likelihood value or the KS statistic. The time series plots of iterations and fitted pdf to the histogram for returns of CAC shares are depicted in Figures \ref{fig3}.
We emphasize that both of log-likelihood and KS statistics reporets in Table \ref{tab1} are evaluated using $\mathsf{STABLE}$. Also the fitted PDF to the histograms in Figure \ref{fig3} are drawn using $\mathsf{STABLE}$ when the PDF parameters are estimated through the EM algorithm. 
\section{Conclusion}
We derive an identity for $\alpha$-stable random variable that is scale-location normal mixture representation. Based on this representation, we propose some EM algorithm to estimate the parameters of $\alpha$-stable distribution. The proposed EM algorithm works very good for admissible ranges of the parameters, i.e., $0<\alpha\leq 2$, $|\beta|\leq 1$, $\sigma>0$, and $\mu_{0} \in \mathbb{R}$. The steps of the proposed EM algorithm are: expectation, maximization, conditionally maximization, and profile log-likelihood maximization. Simulation studies reveal that the EM-based estimators of the tail thickness, skewness, and scale parameters outperform the corresponding maximum likelihood (ML) estimators when scale parameter is large. Also, the EM-based estimator of the scale parameter outperforms the corresponding ML estimator when scale parameter is small. This is not surprising since $\mathsf{STABLE}$ computes the evaluated ML estimators not the exact ML estimators. The sample quantile and empirical characteristic function approaches were eliminated by competitions since the EM approach outperforms them. The simulations reveal that proposed EM algorithm is robust with respect to the initial values. Three sets of real data are used to demonstrate that the proposed EM estimators are close to that of the ML in the sense of goodness-of-fit measures including the log-likelihood and Kolmogorov-Smirnov statistics. A great advantage of the proposed EM algorithm over the ML method is that the proposed EM algorithm can be applied to estimate the parameters of mixture of $\alpha$-stable distributions. Further, since representation (\ref{asymm}) can be adopted for the multivariate case, the proposed EM can be applied for the multivariate $\alpha$-stable distribution and hence is an appropriate approach for modelling the returns of the dependent assets which have multivariate $\alpha$-stable distribution. Two above privileges of the proposed EM algorithm can be considered as a possible future work.
\newpage{}
\noindent{}
{\bf{\Large{Appendix 1: Proof of Proposition 1.1}}}\\
\\
Let $U\sim S_{1}(\alpha,0,1,0)$ and $V\sim S_{1}(\alpha,1,1,0)$ denote two independent $\alpha$-stable random variables. Define $Y=\sigma \bigl(1-|\beta|\bigr)^{\frac{1}{\alpha}}U+\sigma\mathrm{sgn}(\beta)|\beta|^{\frac{1}{\alpha}} V+\mu_{0}-\beta \sigma \tan\bigl(\frac{\pi \alpha}{2}\bigr)$. We have:
\begin{align}
{\varphi}_{Y}{(t)}=&E \exp \biggl\{jt\bigg[\sigma\left(1-\beta\right)^{{\frac{1}{\alpha}}}U+\sigma\mathrm{sgn}(\beta)|\beta|^{{\frac{1}{\alpha}}}V+\mu_{0} -\sigma \beta \tan\Bigl(\frac{\pi \alpha}{2}\Bigr)\bigg]\biggr\} \nonumber \\
=&E \exp \Bigl\{jt\sigma\left(1-|\beta|\right)^{{\frac{1}{\alpha}}}U\Bigr\} E \exp \biggl\{jt\sigma\mathrm{sgn}(\beta)|\beta|^{\frac{1}{\alpha}}V+j t\mu_{0} -j t \beta \sigma \tan\Bigl(\frac{\pi \alpha}{2}\Bigr)\biggr\} \nonumber \\
=&\exp\biggl\{-(1-|\beta|)\left| \sigma t \right|^\alpha-|\beta|\left| \sigma t \right|^\alpha\biggl[1-j\mathrm{sgn}(t\beta)\tan\left(\frac{\pi \alpha}{2}\right)\biggr]+j t\mu_{0}-j t \beta \sigma \tan\Bigl(\frac{\pi \alpha}{2}\Bigr)\biggr\}
\nonumber \\
=&\exp\biggl\{ -\left| \sigma t \right|^\alpha\biggl[1-j\mathrm{sgn}(t)\beta\tan\left(\frac{\pi \alpha}{2}\right)\biggr]+j t\mu_{0}-j t \beta \sigma \tan\Bigl(\frac{\pi \alpha}{2}\Bigr)\biggr\},\nonumber 
\end{align}
where in above, to arrive at the expression after third equality, we use this fact that if $V\sim S_{1}(\alpha,1,1,0)$, then $\sigma\mathrm{sgn}(\beta)|\beta|^{\frac{1}{\alpha}}V\sim S_{1}(\alpha,\mathrm{sgn}(\beta),\sigma |\beta|^{\frac{1}{\alpha}},0)$. The expression after the last equality is the chf of $S_{0}(\alpha,\beta,\sigma,\mu_{0})$. Now, with taking account into the fact that if $U\sim S_{1}(\alpha,0,1,0)$, then $U$ can be represented as a Gaussian scale mixture model, i.e., $U \sim {\cal{N}}(0,2P)$ where $P\sim S_{1}\bigl(\frac{\alpha}{2},1,\bigl(\cos(\pi \alpha/4)\bigr)^{2/\alpha},0\bigr)$, see \cite[p. 20]{Samorodnitsky1994}. Finally, set
$\eta=\sigma\left(1-|\beta|\right)^{\frac{1}{\alpha}}$, $\theta=\sigma \mathrm{sgn}(\beta)|\beta|^{\frac{1}{\alpha}}$, and $\lambda=\beta \sigma \tan\bigl(\frac{\pi \alpha}{2}\bigr)$ to see that $Y\mathop=\limits^d\eta\sqrt{2P}N+ \theta V+\mu_{0}- \lambda$. The proof is complete. 
\\\\\\
{\bf{\Large{Appendix 2: Proof of Proposition 1.2}}}\\
\\
Suppose $E\sim {\cal{E}}(1)$ and $P \sim S_{1}\bigl(\frac{\alpha}{2},1,\bigl(\cos(\pi \alpha/4)\bigr)^{{2/\alpha}},0\bigr)$. Define $R=\frac{E}{P}$, to see that
\begin{align}
P(R\leq r)=&\int_{0}^{\infty}P(E\leq rp) f_{P}(p)dp
=1-\int_{0}^{\infty} \exp\{-rp\}f_{P}(p)dp, \nonumber 
\end{align}
where the last integral in above, i.e., the Laplace transform of random variable $P$ is $\exp\bigl\{-r^{\frac{\alpha}{2}}\bigr\}$; for $r>0$, see \cite[pp. 15]{Samorodnitsky1994}. This means that
$R\sim {\cal{W}}\bigl(\frac{\alpha}{2},1\bigr)$ and consequently
\begin{align}\label{ratio}
\frac{1}{\sqrt{R}}\mathop=\limits^d \frac{1}{W}, 
\end{align}
where $W\sim {\cal{W}}(\alpha,1)$. On the other hand, as noted in Appendix 1, we can write $S=\sqrt{2P}N\sim S_{1}(\alpha,0,1,0)$ where $P\sim S_{1}\bigl(\alpha/2,1,\bigl(\cos(\pi \alpha/4)\bigr)^{2/\alpha},0\bigr)$ and $N \sim {\cal{N}}(0,1)$ are independent. It follows from (\ref{ratio}) that
\begin{align}
\frac{S}{ \sqrt{2E}}\mathop=\limits^d \frac{N}{W}.\nonumber 
\end{align}
The result follows.\\\\\\
{\bf{\Large{Appendix 3: Proof of Proposition 1.3}}}\\
\\
As the main part of the CM-step, implementing the SEM algorithm requires to simulate from posterior distribution of $W_i$ given ${y_i^{**}}$ and $\alpha^{(t)}$. For this purpose, we use the Metropolis-Hasting algorithm. As the candidate, we use the Weibull distribution with the shape parameter $\alpha^{(t)}$. Hence, the acceptance rate $A_{w_i}$, becomes 
\begin{eqnarray} 
A_{w_i}=\min\left\{1, \frac{w^{\text{new}}_i\exp\left\{-\frac{(y_i^{**}w^{\text{new}}_i)^2}{2}\right\}} 
{w^{(t)}_i\exp\left\{-\frac{(y_i^{**}w^{(t)}_i)^2}{2}\right\}}\right\}.\nonumber 
\end{eqnarray} 
Employing the Metropolis-Hasting algorithm, $w^{\text{new}}$ has a little chance for acceptance in each iteration when $|y^{**}_i|$ gets large. Therefore, we use a rejection-acceptance sampling scheme in $m$th cycle of the SEM algorithm at the $t$-th iteration of the proposed EM algorithm to generate from the posterior distribution of $W_i$ given $y^{**}_{i}$ and $\alpha^{(t)}$; for $i=1,\dots,n$. For this, we notice that $f_{W_i|Y^{**}_i}\bigl(w_i\big|y^{**}_i,\alpha^{(t)}\bigr)\propto f_{Y^{**}_i|W_i}(y^{**}_i|w_i) f_{W_i}(w_i)$ and the pdf $f_{Y^{**}_i|W_i}(y_i^{**}|w_i)$ 
is bounded by some value independent of $w_i$, i.e., 
\begin{eqnarray} 
f_{Y^{**}_i|W_i}(y^{**}_i|w_i) \leq \frac{\exp\{-0.5\}}{\sqrt{2\pi}|y^{**}_i|}.\nonumber 
\end{eqnarray} 
So, the rejection-acceptance sampling scheme to generate from $f_{W_i|Y^{**}_i}\bigl(w_i\big|y^{**}_i,\alpha^{(t)}\bigr)$ is given by the following.
\begin{enumerate}
\item Generate a sample from ${\cal{W}}\bigl(\alpha^{(t)},1\bigr)$, say $w_i$.
\item Generate a sample from uniform distribution over $\left(0,\exp\{-0.5\}/(\sqrt{2\pi}|y^{**}_{i}|)\right)$, say $u$.
\item If $u<\frac{w_i}{\sqrt{2\pi}}\exp\left\{-\frac{{y^{**}_{i}}^2 w_{i}^{2}}{2}\right\}$, then accept $w_i$ as a realization of $f_{W_i|Y^{**}_i}(w_i|y^{**}, \alpha^{(t)})$; otherwise, start from step 1. 
\end{enumerate}
{\bf{\Large{Appendix 3: Evaluating $E^{(t)}_{ri}$}}}\\
\\
At $t$-th iteration of the proposed EM algorithm, for $i$th observed value $y_{i}$, define three $K\times K$ matrices such as A, B, and C. Elements of matrix A are independent realizations from $S_{1}\bigl(\alpha^{(t)}/2,1,\bigl(\cos(\pi \alpha^{(t)}/4)\bigr)^{2/\alpha^{(t)}},0\bigr)$ and elements of matrix B are coming from $S_{1}(\alpha^{(t)},1,1,0)$. 
Assume that $A_{rc}$ and $B_{rc}$, respectively, are the $c$th element of the $r$th row of matrices A and B, then the $c$th element of the $r$th row of matrix C, i.e., $C_{rc}$ is $\phi\bigl(y_{i},\mu^{(t)}_{0}-\lambda^{(t)}+\eta^{(t)}B_{rc}, 2A_{rc}(\theta^{(t)})^2\bigr)$. Now, 
\begin{align}
E^{(t)}_{0i}=&\sum_{r=1}^{K}\sum_{c=1}^{K}\frac{C_{rc}}{B_{rc}K^2},\nonumber\\
E^{(t)}_{1i}=&\sum_{r=1}^{K}\sum_{c=1}^{K}\frac{A_{rc}C_{rc}}{B_{rc}K^2},\nonumber\\
E^{(t)}_{2i}=&\sum_{r=1}^{K}\sum_{c=1}^{K}\frac{A^{2}_{rc}C_{rc}}{B_{rc}K^2}.\nonumber
\end{align}
Also, $f_{Y}\bigl(y_i\big|{\alpha}^{(t)},\beta^{(t)}, {\sigma}^{(t)}, {\mu_{0}}^{(t)}\bigr)$ is approximated as
\begin{align}
f_{Y}\bigl(y_i\big|{\alpha}^{(t)},\beta^{(t)}, {\sigma}^{(t)}, {\mu_{0}}^{(t)}\bigr) \approx&
\sum_{r=1}^{K}\sum_{c=1}^{K}\frac{C_{rc}}{K^2}.\nonumber
\end{align} 
It should be noted that the constant K must be large enough. Using $K=100$, the quantities $E^{(t)}_{ri}$; for $r=0,1,2$ and $i=1,\dots,n$, are approximated very accurately.

\end{document}